\documentclass{amsart}
\usepackage[curve,matrix,arrow]{xy}
\usepackage{amssymb, amscd}
\newtheorem{theorem}{Theorem}[section]
\newtheorem{lemma}[theorem]{Lemma}
\newtheorem{corollary}[theorem]{Corollary}
\newtheorem{proposition}[theorem]{Proposition}
\newtheorem{definition}[theorem]{Definition}
\newtheorem{remark}[theorem]{Remark}

\newtheorem{assumptions}[theorem]{Standing assumptions}
\newtheorem{numbered equation}[theorem]{(\hspace*{-0.1cm}}

\DeclareMathOperator{\Hom}{Hom}

\DeclareMathOperator{\map}{map}

\DeclareMathOperator{\Ho}{Ho}
\DeclareMathOperator{\HoC}{Ho(\C)}
\DeclareMathOperator{\hocolim}{hocolim}
\DeclareMathOperator{\colim}{colim}
\DeclareMathOperator{\Ev}{Ev_0}

\DeclareMathOperator{\Sym}{Sym}
\DeclareMathOperator{\co}{C_0}

\newcommand{\pPsi}{\Theta}

\newcommand{\mIC}{{\mathbb I}_\C}
\newcommand{\mID}{{\mathbb I}_\D}

\newcommand{\bS}{{\mathbb S}}

\newcommand{\mQ}{{\mathbb Q}}
\newcommand{\bQ}{{\mathbb Q}}
\newcommand{\bQt}{\widetilde{\mathbb Q}}
\newcommand{\bI}{{\mathbb I}}
\newcommand{\hh}{{\mathbb H}}
\newcommand{\mZ}{{\mathbb Z}}
\newcommand{\bZ}{{\mathbb Z}}
\newcommand{\bZt}{\widetilde{\mathbb Z}}
\newcommand{\tZ}{\widetilde{\mathbb Z}}
\newcommand{\tU}{\widetilde{U}}

\newcommand{\A}{\mathcal {A}}
\newcommand{\N}{\mathcal {N}}
\newcommand{\cA}{{\mathcal A}}

\newcommand{\cR}{{\mathcal R}}
\newcommand{\C}{{\mathcal C}}
\newcommand{\cC}{{\mathcal C}}
\newcommand{\D}{{\mathcal D}}
\newcommand{\E}{\mathcal E}

\newcommand{\G}{\mathcal G}

\newcommand{\simp}{{\mathcal S}_*}
\newcommand{\simpu}{\mathcal S}
\newcommand{\J}{\mathcal F}

\newcommand{\alg}{{\mbox{-}\mathcal A}lg}

\newcommand{\spec}{Sp^{\Sigma}}
\newcommand{\Spec}{Sp^{\Sigma}}

\newcommand{\cD}{\D}

\newcommand{\sset}{\mathcal{S}_*}

\newcommand{\boxprod}{\mathbin{\square }}

\newcommand{\iso}{\cong}

\newcommand{\sm}{\wedge}
\newcommand{\Smash}{\wedge}

\newcommand{\tensor}{\otimes}
\newcommand{\tens}{\otimes}

\newcommand{\widebar}{\overline}
\newcommand{\mc}{\colon \,}

\newcommand{\proj}{\mbox{\scriptsize{proj}}}
\newcommand{\inj}{\mbox{\scriptsize{inj}}}
\newcommand{\id}{\mbox{Id}}
\newcommand{\Mod}{\mbox{-}{\mathcal M}od}
\newcommand{\Alg}{\mbox{-}{\mathcal A}lg}
\newcommand{\monoid}{\mbox{-}{\mathcal M}onoids}
\newcommand{\mon}{\mbox{\scriptsize{mon}}}
\newcommand{\Modr}{{\mathcal M}od\mbox{-}}
\newcommand{\sab}{s{\mathcal A}b}
\newcommand{\ch}{ch^+}
\newcommand{\chZ}{{\mathcal C}h}
\newcommand{\Ch}{{\mathcal C}h}
\newcommand{\dga}{{\mathcal DGA}}

\newcommand{\dgar}{{\mathcal DGA}_{R}}

\newcommand{\DGAZ}{{\mathcal DGA}_{\mZ}}
\newcommand{\Gam}{\Gamma}

\renewcommand{\to}{\longrightarrow}
\newcommand{\ra}{\longrightarrow}
\newcommand{\lla}{\longleftarrow}
\newcommand{\cofib}{\rightarrowtail}
\newcommand{\varrow}[1]{\hbox to #1{\rightarrowfill}}
\newcommand{\varl}[2]{\stackrel{#2}{\hbox to #1{\leftarrowfill}}}

\newcommand{\varrx}[2]{\stackrel{#2}{\hbox to #1{\rightarrowfill}}}
\newcommand{\parallelarrows}[1]{\begin{array}{c} {\hbox to
#1{\rightarrowfill}}  \vspace{-0.35cm} \\ {\hbox to
#1{\rightarrowfill}} \end{array}}
\newcommand{\adjoint}{\rightleftarrows}

\begin{document}

\title[DGAs and $H{\mathbb Z}$-algebras]{ 
$H{\mZ}$-algebra spectra are 
differential graded algebras 
}

\date{\today; 2000 AMS Math.\ Subj.\ Class.: 55P43, 18G35, 55P62, 18D10} 
\author{Brooke Shipley}
\thanks{
Research partially supported by NSF Grants No. 0417206 and No. 0134938 and a
Sloan Research Fellowship}
\address{Department of Mathematics \\ 510 SEO m/c 249\\   
851 S. Morgan Street\\ Chicago, IL 60607-7045 \\ USA}
\email{bshipley@math.uic.edu}
\maketitle

%%%%%
{\small Abstract:  
We show that the homotopy theory of differential graded algebras
coincides with the homotopy theory of $H\mZ$-algebra spectra.
Namely, we construct Quillen equivalences between the Quillen model categories 
of (unbounded) differential graded algebras and $H\mZ$-algebra spectra.
We also construct Quillen equivalences between the 
differential graded modules and module spectra over these algebras.
We use these equivalences in turn to produce algebraic models for rational 
stable model categories.  We show that basically 
any rational stable model category is Quillen equivalent to modules 
over a differential graded $\bQ$-algebra (with many objects).
}

\section{Introduction}\label{sec-intro}

In this paper we show that the algebra of spectra arising in 
stable homotopy theory encompasses homological algebra.  
Or, from another perspective, this shows that the algebra of spectra is a
generalization of homological algebra.
In the algebra of spectra, spectra take the place of abelian groups 
and the analogue of the tensor product is a {\em symmetric monoidal 
smash product}.  Using this smash product, the definitions
of module spectra, associative and commutative ring spectra, and algebra
spectra are easy categorical exercises.  The main difference here is 
that in the algebra of spectra the sphere spectrum, $\bS$, which encodes all of 
the information about the stable homotopy groups of spheres, takes the place 
of the integers, the free abelian group on one generator.

Originally, 
ring spectra (or $A_{\infty}$-ring spectra) were built from topological spaces 
with ring structures which are associative up to coherent homotopies. 
Since a discrete ring is strictly associative, rings are examples of ring 
spectra; associated to any ring $R$ is the Eilenberg-Mac Lane ring 
spectrum $HR$.  One good expository introduction to spectra and ring spectra
appears in~\cite{green-expos}. 
As in classical
algebra, to understand ring spectra we study their modules.  
Since we are considering rings up to homotopy though, we focus 
on the type of homological invariants and derived functors best captured 
by the derived category rather than the category of modules itself.  
Robinson, in~\cite{robinson}, defined a notion of $A_{\infty}$-modules 
over an Eilenberg-Mac Lane spectrum $HR$ and a notion of
homotopy.  He then showed that up to homotopy this category of
$A_{\infty}$-modules is equivalent to $\cD(R)$, the unbounded derived category
%or the homotopy category of (unbounded) differential graded $R$-modules.
of $R$.  The complicated notions of algebraic structures up to coherent
homotopies prevented anyone from generalizing this result to algebras 
or even defining a notion of $A_{\infty}$-$HR$-algebra spectra. 

In the modern settings of spectra referred to in the first paragraph, 
defining algebra spectra is simple.  Here the sphere spectrum is
the initial ring and spectra are $\bS$-modules, the modules over the 
sphere spectrum $\bS$.  Ring spectra are then the (strictly associative) 
$\bS$-algebras.  For a commutative 
$\bS$-algebra $A$, $A$-algebra spectra are the monoids in the category of 
$A$-modules.       %~\cite{ekmm, hss}.
Given a discrete ring $R$, one can construct $HR$ as an $\bS$-algebra which 
is commutative if $R$ is~\cite[1.2.5]{hss}. 
Since Robinson showed that $HR$-module spectra agree with differential
graded $R$-modules ``up to homotopy," one might guess that $HR$-algebra 
spectra should capture the same ``up to homotopy" information as differential 
graded $R$-algebras.  Verifying this conjecture is the subject of this paper. 

To formulate this statement rigorously we use Quillen model category
structures.
One may define the derived category of a ring $R$ as the category of
dg $R$-modules with quasi-isomorphisms formally inverted.  The problem with 
this formal process of inversion is that one loses control of the morphisms; 
the morphisms between two objects may not form a set.  Quillen model structures
give a way around this problem.  A {\em Quillen model structure}
on a category $\cC$ is a choice of three subcategories called weak equivalences,
cofibrations and fibrations which satisfy certain axioms~\cite{Q}.  
\cite{ds} is a good introduction to model categories; our standard
reference though is~\cite{hovey-book}.  Given a
Quillen model structure, inverting the weak equivalences, $W,$ produces a 
well-defined homotopy category $\Ho(\cC) = \cC[W^{-1}]$.   
%The categories of module spectra and algebra
%spectra mentioned above all have model structures.  
The category of dg $R$-modules, $R\Mod$, forms a model category and one 
can show that $\Ho(R\Mod)$ is equivalent to
$\cD(R)$.   A {\em Quillen adjoint pair} (or Quillen map) between two Quillen 
model categories $\cC$ and
$\cD$ is an adjoint pair of functors $F\mc \cC \rightleftarrows \cD\mc G$ which
preserves the model structures and thus induces adjoint derived functors on the
homotopy categories $LF \mc \Ho(\C) \rightleftarrows \Ho(\cD) \mc RG$.  A 
{\em Quillen equivalence} is a Quillen map such that the derived functors  
induce an equivalence on the homotopy categories.  Two model categories are 
{\em Quillen equivalent} if they can be connected by a zig-zag of Quillen 
equivalences.  
Quillen equivalent model categories represent the same `homotopy theory.'
A subtle but important point is that a Quillen equivalence is 
stronger than just having an equivalence of the associated homotopy 
categories; see examples at~\cite[3.2.1]{stable} or~\cite[4.5]{morita}.  

In~\cite[5.1.6]{stable}, we strengthened Robinson's result mentioned above by 
showing that the model category of $HR$-module spectra is Quillen equivalent to
$R\Mod$.   It follows that $\cD(R)$, $\Ho(R\Mod)$ and $\Ho(HR\Mod)$
are all equivalent.  Connective versions of
these results for modules also appeared in~\cite[4.4, 4.5]{schwede-thesis}; 
see also~\cite[1.1]{emmc}.  \cite[IV.2.4]{ekmm} shows that $\Ho(HR\Mod)$
and $\cD(R)$ are equivalent categories although a Quillen equivalence is 
not established.  A similar homotopy category level equivalence is 
stated in~\cite{stanley} between the associated algebra categories
$\Ho(HR\Alg)$ and $\Ho(\dgar)$. 

If $R$ is a commutative ring, 
$R\Mod$ and $HR\Mod$ both have symmetric monoidal products, 
but the Quillen maps used in~\cite{stable} do not respect these products  
and hence do not induce functors on the associated 
categories of algebras.  Another approach is thus needed if one wants to 
compare the homotopy theory of $HR$-algebra spectra, $HR\Alg$, and 
dg $R$-algebras, $\dgar$.  The main result of this paper is then the 
construction of a new zig-zag of Quillen equivalences between $R\Mod$ and 
$HR\Mod$ which preserve the relevant product structures and induce Quillen 
equivalences between $HR\alg$ and $\dgar$. 

\begin{theorem}\label{thm-ch-hZ}      
For any discrete commutative ring $R$, 
the model categories of 
(unbounded) differential graded
$R$-algebras and $HR$-algebra spectra are Quillen equivalent.
The associated composite derived functors are denoted
$\hh\mc \dgar \to HR\alg$ and $\Theta\mc
HR\alg \to \dgar$.
\end{theorem}       

The proof of this theorem appears in Section~\ref{2}.
We also construct Quillen equivalences between the categories of modules over
these algebras; see Corollary~\ref{Big Cor}.
We then use these Quillen equivalences 
to construct algebraic models for {\em rational stable model categories};
see Corollary~\ref{cor-rational}.   

In~\cite{man-taq}, an analogue of Theorem~\ref{thm-ch-hZ} is proved for 
$E_{\infty}$-algebras, the commutative analogue of $A_{\infty}$-algebras
which are associative and commutative up to coherent homotopies.  Mandell 
shows that the homotopy category of $E_{\infty}$-$HR$ algebras is equivalent 
to the homotopy category of dg $E_{\infty}$-$R$-algebras~\cite[7.11]{man-taq}.  
He also claims these equivalences
can be lifted to Quillen equivalences using the techniques of~\cite{emmc}.  
The model categories of $E_{\infty}$-$HR$ algebras and commutative 
$HR$-algebras are shown to be Quillen equivalent in ~\cite[7.2]{man-taq}.   
If $R$ is not an algebra over the rationals though, the category of 
commutative dg $R$-algebras does not even have a well behaved homotopy theory;
that is, there is no model category on commutative dg $R$-algebras with weak 
equivalences and fibrations determined on the underlying
dg $R$-modules (the quasi-isomorphisms and surjections).  

In general then, the Quillen equivalence in Theorem~\ref{thm-ch-hZ} cannot be  
extended to categories of commutative algebras.   Rationally though, it
should hold.  We do not consider this extension here though because it would 
require different techniques.  The following related statement,
however, is simple to prove here and is used in the construction of an
algebraic model of rational $T^n$-equivariant spectra in~\cite{gs}.  Its
proof appears near the end of Section~\ref{sec-proof}.
Recall $\pPsi$ from Theorem~\ref{thm-ch-hZ}. 

\begin{theorem}\label{thm-Q-comm}
For $C$ any commutative $H\bQ$-algebra, $\pPsi C$ is naturally weakly equivalent
to a commutative differential graded $\bQ$-algebra.
\end{theorem}
 
{\em Acknowledgments:} I would like to thank Stefan Schwede for
many fruitful conversations on this material.  
In fact, this project grew out of a sequence of
projects with him.  In particular, Corollary~\ref{cor-rational}
is really unfinished business from~\cite{stable} and even
appeared in various preprint versions.  I would also like to thank 
Dan Dugger and Mike Mandell for many helpful conversations.  
I would like to thank John Greenlees, David Gepner and the 
referee for their careful reading of an earlier version of this paper. 

\section{$\mZ$-graded chain complexes and Eilenberg-Mac Lane spectra}\label{2}

Throughout this paper let $R$ be a discrete commutative ring.
In this section we define the functors involved in comparing differential
graded $R$-algebras (DGAs) and $HR$-algebra spectra.   As mentioned in
the introduction, the functors used in previous comparisons of
dg $R$-modules and $HR$-module spectra do not induce functors on the
associated categories of algebras because they do not respect the product 
structures on these module categories.  The main input for showing 
the categories of algebras are Quillen equivalent is then a new 
comparison of the module categories via functors which do respect
the product structures.   Namely, we construct a zig-zag of Quillen
maps between $R\Mod$ and $HR\Mod$ in which the right adjoints are
lax monoidal functors.  More precisely, we show that each of the
adjoint pairs in this zig-zag are {\em weak monoidal Quillen equivalences} in 
the sense of~\cite{emmc}. The results of~\cite{emmc} then show that they
induce the required Quillen equivalences of categories of algebras and modules. 

Here a {\em monoid} (or algebra) is an object $R$ together with a
``multiplication'' map $R\tensor R \xrightarrow{\mu} R$ and a 
``unit'' $\mathbb I\xrightarrow{\eta} R$ which satisfy
certain associativity and unit conditions (see [MacL, VII.3]).
If $R$ is a monoid, a {\em left $R$-module} (``object with left $R$-action'' 
in [MacL, VII.4]) is an object $N$ together with an action map 
$R\tensor N\xrightarrow{\alpha} N$ satisfying
associativity and unit conditions (see again [MacL, VII.4]).

We begin this section by briefly recalling the necessary definitions and 
results from~\cite{emmc}.   After this we define the categories and
functors involved in the zig-zag of Quillen maps between $R\Mod$ and
$HR\Mod$.   
We defer the verification that these functors are weak monoidal 
Quillen pairs to Section~\ref{sec-proof}.  
We then prove our main results about the Quillen equivalence of categories 
of algebras and modules arising in spectral 
algebra and homological algebra.   We end this section by using these 
equivalences to construct algebraic models for rational stable model
categories. 

\subsection{Weak monoidal Quillen pairs}
Although our main focus here is describing conditions under
which a Quillen map induces a Quillen equivalence on the associated
categories of algebras and modules, we first need to consider conditions
on the categories themselves.  Namely, the categories of algebras and
modules must have model structures created from the model structure
on the underlying category.    In~\cite[4.1]{ss} it is shown that
if every object in $\cC$ is small and $\cC$ is a cofibrantly generated, 
monoidal model category which satisfies the monoid axiom, then the categories 
of monoids, modules and algebras over $\cC$ have model structures created
by the forgetful functor to $\cC$, that is, with
the weak equivalences and fibrations determined on $\cC$.   We now
recall the definitions of a monoidal model category and the monoid axiom.
We refer
to~\cite[2.2]{ss} or~\cite[2.1.17]{hovey-book} for the definition of
cofibrantly generated model categories.  In Section~\ref{sec-monoid},
we verify these conditions for the categories considered in this paper. 

The requirement that $\cC$ is a {\em monoidal model category} ensures  
that $\cC$ is a monoidal category where the monoidal product is compatible with 
the model
structure.  In general, this would also require a compatibility condition
for the unit~\cite[4.2.6 (2)]{hovey-book}, but this condition is not necessary when the unit 
is cofibrant which is the case for all of the categories considered in this 
paper.   

\begin{definition}\label{def-mmc}\cite[3.1]{ss} {\em
A model category $\C$ is a {\em monoidal model category} if it has a closed 
symmetric monoidal structure with product $\otimes$ and cofibrant unit object 
$\bI_{\C}$ and satisfies the following axiom.

\noindent
{\em Pushout product axiom.} Let $f\mc A \to B$ and $g\mc K \to L$ be 
cofibrations in $\C$. Then the map 
$$f\boxprod g \mc A \otimes L   \coprod_{  A\otimes K}    B\otimes   K \to 
B \otimes L$$
is also a cofibration. If in addition one of the former maps is a weak 
equivalence, so is the latter map.}
\end{definition}

The {\em monoid axiom} is the extra condition used in~\cite{ss} to 
extend a model structure on a monoidal model
category to model structures on the associated categories of algebras
and modules. 

\begin{definition}\label{def-ma}\cite[3.3]{ss} 
{\em A monoidal model category $\C$ satisfies the {\em monoid axiom} if
any map obtained by pushouts and (possibly transfinite) compositions
from maps of the form $f \sm \id_Z: A \sm Z \to B \sm Z$ for $f$
a trivial cofibration and $Z$ any object in $\C$ is a weak equivalence. 
}\end{definition}

Now we return to the conditions on a Quillen map necessary for
inducing a Quillen equivalence on the associated categories of
algebras and modules.   The main required property is that the
right adjoint is lax monoidal and induces a Quillen equivalence on
the underlying categories.   Additional homotopy properties on the
units and the left adjoint are also required though.  

\begin{definition}\label{def-wmQp}\cite{emmc}
\em{
A Quillen adjoint pair $(\lambda, R)$ between
monoidal model categories $(\C,\otimes)$ and $(\D, \sm)$ with
cofibrant unit objects, $\bI_{\C}$ and $\bI_{\D}$, is a 
{\em weak monoidal Quillen
pair} if the right adjoint $R$ is lax monoidal~\cite[20.1]{mmss}, 
$\tilde{\nu} \mc
\lambda(\bI_{\D}) \to \bI_{\C}$ is a weak equivalence and $\tilde{\varphi}\mc
\lambda(A \sm B) \to \lambda A \otimes \lambda B$ is a weak equivalence
for all cofibrant objects $A$ and $B$ in $\D$.  See~\cite[3.6]{emmc}
for the definition of $\tilde{\nu}$ and $\tilde{\varphi}$, the co-monoidal 
structure on $\lambda$ adjoint to the monoidal structure on $R$. 
A {\em strong monoidal Quillen pair} is a weak pair for which 
$\tilde{\nu}$ and $\tilde{\varphi}$ are isomorphisms; for example, this
holds when  $\lambda$ is strong monoidal~\cite[20.1]{mmss}.  
When $(\lambda, R)$ is a Quillen
equivalence, we change these names accordingly.}
\end{definition}

In~\cite[3.3]{emmc} it is shown that if $\lambda \mc \cD \adjoint \cC\mc R$ is 
a weak monoidal Quillen pair then the monoid valued lift
$R\mc \cC\monoid \to \cD\monoid$ has a left adjoint $L^{\mon} \mc
\cD\monoid \to \cC\monoid$.   $L^{\mon}$ does not agree with
$\lambda$ on underlying objects in general unless $(\lambda, R)$ is a strong 
monoidal Quillen pair.   In the strong monoidal case we usually abuse notation
and denote the induced functors on monoids as $\lambda$ and $R$ again. 

\begin{assumptions}\label{stass}{\em
In the following theorem, assume that in $\cC$ and $\cD$  
every object is small and $\cC$ and $\cD$ are cofibrantly generated, 
monoidal model categories which satisfy the monoid axiom.  Assume
further that the unit objects in $\cC$ and $\cD$ are cofibrant.  We will
refer to these as our {\em standing assumptions}.  It follows 
by~\cite[4.1]{ss} that model structures on monoids in 
$\cC$ and $\cD$ are created by the forgetful functors. 
}\end{assumptions}

\begin{theorem}\cite[3.12]{emmc}\label{thm-emmc} 
Let $\lambda \mc \cD \to \cC\mc R$ be a weak monoidal Quillen
equivalence.  If $\cC$ and $\cD$ satisfy the standing assumptions
listed above, 
then the adjoint functor pair
\[ L^{\mon}\mc \cD\monoid \adjoint \cC \monoid \mc R \]
is a Quillen equivalence between the respective model categories of
monoids.   If $R$ is the right adjoint of a strong monoidal Quillen
equivalence, then $L^{\mon}$ agrees with $\lambda$ on underlying
objects.   
\end{theorem}
    
\begin{remark}
{\em 
Unfortunately, we have not found a single weak monoidal Quillen equivalence 
between $HR$-module spectra and differential graded $R$-modules.
Instead, we construct a three step zig-zag of weak monoidal Quillen 
equivalences; see the subsection below. 
For strong monoidal Quillen pairs, the functors induce adjoint pairs on the 
associated categories of modules and algebras and it is relatively easy to 
show these are Quillen equivalences, following~\cite[Section 16]{mmss} for 
example.  Two of our three steps will be via strong monoidal Quillen
equivalences.  For the middle step though, the argument is a bit more 
involved.  Thus, it turns out to be simplest for all three steps to use the 
general framework of weak monoidal Quillen equivalences~\cite[3.12]{emmc} for 
extending Quillen equivalences to categories of algebras and modules.
}\end{remark}

\subsection{The zig-zag of functors}
Here we describe the categories and functors involved in the zig-zag
of Quillen maps between $HR\Mod$ and $R\Mod$.  In Section~\ref{sec-monoid}
we show that the Standing assumptions~\ref{stass} apply to
these categories and in Section~\ref{sec-proof} we show that each of the
adjoint functor pairs is a weak monoidal Quillen equivalence. 
For simplicity we concentrate on $R = \bZ$ but in every case $\bZ$ 
could be replaced by any discrete commutative ring.  We denote
the category $\bZ\Mod$ by $\chZ$.  
First we describe the two intermediate categories which are both analogues of
symmetric spectra as considered in~\cite[Section 7]{hovey-stable}.  We   
recall the basic definitions here.  Let $\C$ be a closed symmetric
monoidal category with monoidal product $\tensor$ and unit $\bI_{\C}$.

\begin{definition}
{\em 
Let $\Sigma$ be the category with objects the sets
$\widebar{n}= \{ 1,2, \dots, n\}$ for $n \geq 0$, (note that 
$\widebar{0}$ is the empty set)
and morphisms the isomorphisms.  The category of {\em symmetric
sequences in $\C$} is the functor category $\C^{\Sigma}$. 
$\C^{\Sigma}$ is a symmetric monoidal category; the monoidal
product is defined by 
\[(X \tensor Y)_n = \bigvee_{p+q =n} (\Sigma_n)_{+} \Smash_{\Sigma_p \times
\Sigma_q} (X_p \tensor Y_q). \]
}\end{definition}  

Given an object $K$ in $\C$, define $\Sym(K)$ as the symmetric sequence
$(\bI_{\C}, K, K \otimes K, \dots, K^{\otimes n}, \dots)$ 
where the symmetric group $\Sigma_n$ acts on $K^{\otimes n}$
by permutation.  $\Sym(K)$ is a commutative monoid in $\C^{\Sigma}$; 
we consider the category of $\Sym(K)$-modules.  

\begin{definition}\label{def-spec}  
{\em The category of {\em symmetric spectra over $\cC$},  $\spec(\C,K)$, is the category of
modules over $\Sym(K)$ in $\C^{\Sigma}$. $\spec(\C,K)$ is a symmetric
monoidal category; the monoidal product $X \Smash Y$ is defined as the 
coequalizer of the two maps
\[ X \tensor \Sym(K) \tensor Y \parallelarrows{.5cm} X \tensor Y.\]} 
\end{definition}  

Let $\simp$, $\sab$ and $\ch$ %and $\chZ$ 
denote the categories of simplicial sets, simplicial abelian groups and 
non-negatively graded chain complexes. 
Let $\tZ\mc \simp \to
\sab$ be the functor such that $\tZ(K)_n$ is the free abelian group
on the non-basepoint simplices in $K_n$.  
Through this paper, let $S^1= \Delta[1]/ \partial \Delta[1]$ denote the 
simplicial circle, $S^n = (S^1)^{\sm n}$ and  
$\mZ[n]$ be the chain complex which contains a single
copy of the group $\mZ$ in dimension $n$. 
We consider $\spec(\sab, \tZ S^1)$ and $\spec(\ch,\mZ[1])$.
To ease notation, we refer to these categories as $\spec(\sab)$ and
$\spec(\ch)$ from now on. %and $\spec(\chZ)$ from now on.

$\spec(\sab)$ and $\spec(\ch)$ are the two intermediate categories 
in our comparison between $H\mZ\Mod$ and $\chZ$.
Next we define the adjoint functors involved in the 
three step comparison between $H\mZ\Mod$ and $\chZ$ as displayed below
with left adjoints on top. 

\[\xymatrix@=10mm{ 
H\mZ\Mod \quad 
\ar@<.4ex>^-{Z}[r] &
\quad  \spec(\sab) \quad 
\ar@<.4ex>^-{U}[l]
\ar@<-.4ex>_-{\phi^*N}[r] &
\quad  \spec(\ch) \quad  
\ar@<-.4ex>_-{L}[l]
\ar@<.4ex>^-{D}[r] &
\quad \chZ \quad \ar@<.4ex>^-{R}[l]
}
\]

\noindent
These functors induce functors on the categories of monoids.
The composite derived functors (compare with~\cite[1.3.2]{hovey-book}) 
mentioned in
the introduction in Theorem~\ref{thm-ch-hZ} are then 
\[ \pPsi = D c \phi^* N  Z c     
\mbox{\  and \ } \hh= U  L^{\mon} c R \] 
where $c$ denotes the cofibrant replacement functors 
in each of the model categories of monoids.  Fibrant replacement
functors are not needed because each of the right adjoints
$U, \phi^* N,$ and $R$ preserve all weak equivalences.  
In Section~\ref{sec-proof}
we show that $(Z,U)$ and $(D,R)$ are both strong monoidal Quillen
equivalences.  It follows that the left adjoints induced on
the categories of monoids are just the restrictions of the underlying
functors.  Thus we still denote the induced
left adjoints by $Z$ and $D$.   Since $(L, \phi^*N)$ is only a weak
monoidal Quillen equivalence, the induced left adjoint on monoids here
is $L^{\mon}$ as discussed above  Theorem~\ref{thm-emmc}.  

First consider the functors between $H\bZ\Mod$ and $\spec(\sab)$.
Since $\tZ$ is strong monoidal and we defined $S^n = (S^1)^{\sm n}$, 
then $(\tZ S^1)^{\otimes n} \iso \tZ S^n$.  Therefore the forgetful functor 
from simplicial abelian groups to pointed simplicial sets takes 
$\Sym(\tZ S^1)= \tZ \bS$ to the symmetric spectrum 
\[H\mZ = (\mZ, \tZ S^1, \tZ  S^2, \cdots, \tZ S^n, \cdots ) \]
as defined in~\cite[1.2.5]{hss}.
This forgetful functor induces the functor $U\mc\spec(\sab )
\to H\mZ\Mod$; we define its left adjoint $Z$ in 
the proof of Proposition~\ref{prop-hz-sab}.  

The next step is to compare $\spec(\sab )$ and $\spec(\ch)$.
Let $N\mc \sab \to \ch$ denote the normalization functor and
$\Gamma$ denote its inverse as defined in~\cite[2.1]{emmc}, for example.  
Applying the normalization functor from $\sab$ to $\ch$ to each
level induces a functor from $\spec(\sab)$ to
the category of modules in $(\ch)^{\Sigma}$ over 
$\N= N(\Sym(\tZ S^1))=N\tZ\bS=(N\tZ S^0, N\tZ S^1, N\tZ S^2, \dots)$.  $\N$ is a commutative monoid
since $N$ is a lax symmetric monoidal functor by~\cite[2.6]{emmc}.  
By identifying
$N\tZ S^1$ as $\mZ[1]$, we see there is a ring map 
$\phi\mc\Sym(\mZ[1]) \to \N$ which in degree $n$ is induced by the monoidal
structure on $N$, $(\mZ[1])^{\otimes n} \iso 
(N\tZ S^1)^{\otimes n} \to N(\tZ S^1{}^{\otimes n})$;  
see~\cite[29.7]{may}. 
The composition of $N$ and forgetting along $\phi$
gives a functor $\phi^* N \mc \spec(\sab) \to \spec(\ch)$.
A left adjoint to $\phi^*N$ exists by~\cite[Section 3.3]{emmc};
denote it by $L$.
Note, $L$ is not isomorphic to the composite of underlying left adjoints,
$\Gam \phi_*$, because the adjoint of the identity on $\N$, namely 
$\Sym(\tZ S^1) \to \Gam N \Sym (\tZ S^1)$, is not
a ring map; see~\cite[2.14]{emmc}.

Define a functor $R: \Ch \to \Spec(\ch)$ by setting 
$(RY)_m = C_0 (Y \tensor \bZ[m])$.  
Here $\bZ[m]$ is the chain complex with a single copy of $\bZ$ 
concentrated in degree $m$. $C_0$ is the connective cover; 
it is the right adjoint to the inclusion of 
$\ch$ in $\Ch$.  
One can check that there
are isomorphisms $(RY)_m \xrightarrow{\iso} C_0(\bZ[-1] \otimes (RY)_{m+1})$
with adjoints $\bZ[1] \otimes (RY)_m \to (RY)_{m+1}$ which provide
the module structure over $\Sym(\bZ[1])$. 
We define $D$, the left adjoint of $R$, in the paragraph 
above Proposition~\ref{prop-DR}. 

%In Propositions~\ref{prop-chs} and~\ref{prop-chZ} we show that
%$i$ and its right adjoint $\co$ and $\Ev$ and its left adjoint $F_0$
%form Quillen equivalences and 
%induce Quillen equivalences on the associated categories of algebras
%and modules. 

\subsection{Statements of results}
Here we summarize the properties of the categories and functors
which are verified in Sections~\ref{sec-monoid} and~\ref{sec-proof}.   We then
prove Theorem~\ref{thm-ch-hZ}.  

The following statement is proved as Corollary~\ref{cor-mon-ax}. 

\begin{proposition}\label{prop-mon-ax-one}
$H\bZ\Mod$, $\spec(\sab)$, $\spec(\ch)$ and $\chZ$ satisfy the Standing 
assumptions~\ref{stass}.  It follows by~\cite[4.1]{ss} 
that there are model structures 
on the categories of modules and algebras over these categories
with fibrations and weak equivalences defined on the underlying
model category.
\end{proposition}

\begin{proposition}\label{prop-list}
The following statements are proved as Propositions~\ref{prop-hz-sab},
~\ref{prop-sab-ch} and~\ref{prop-DR2}. 
%and~\ref{prop-chZ}.
\begin{enumerate}
\item $( \spec(\sab), H\mZ\Mod, Z, U)$ is a strong monoidal Quillen equivalence.

\item 
$( \spec(\sab), \spec(\ch), L, \phi^* N)$ is a weak monoidal Quillen equivalence.  

\item $(\Ch, \spec(\ch), D, R)$ is a strong monoidal Quillen equivalence. 
% \item $(\spec(\chZ), \spec(\ch), i, \co)$ is a strong monoidal Quillen equivalence. 

%\item $ (\spec(\chZ),  \chZ, F_0, \Ev)$ is a strong monoidal Quillen equivalence. 
\end{enumerate}
\noindent
Moreover, the right adjoint $(U, \phi^*N, R)$ in each of these pairs 
preserves all weak equivalences. 
\end{proposition}

Theorem~\ref{thm-ch-hZ} then follows from Theorem~\ref{thm-emmc}. 

\begin{proof}[Proof of Theorem~\ref{thm-ch-hZ}]
Proposition~\ref{prop-list} and Proposition~\ref{prop-mon-ax-one}
verify all of the hypotheses required to apply Theorem~\ref{thm-emmc}
to the three step zig-zag between $H\bZ\Mod$ and $\chZ$.  
These three adjoint pairs thus induce Quillen equivalences on
the categories of monoids.    
\[\xymatrix@=10mm{ 
H\mZ\Alg \quad 
\ar@<.4ex>^-{Z}[r] &
\quad  \spec(\sab)\monoid \quad 
\ar@<.4ex>^-{U}[l]
\ar@<-.4ex>_-{\phi^*N}[r] &
\quad  \spec(\ch)\monoid \quad  
\ar@<-.4ex>_-{L^{\mon}}[l]
\ar@<.4ex>^-{D}[r] &
\quad \DGAZ \quad \ar@<.4ex>^-{R}[l]
}
\]
Note, monoids in $H\bZ\Mod$ are the $H\mZ$-algebra spectra, denoted $H\mZ\Alg$
above.   Also, the monoids in $\chZ$ are differential graded algebras, $\DGAZ$. 
For an arbitrary discrete commutative ring $R$ one can replace $\mZ$ and 
abelian groups by $R$ and $R$-modules in all of the statements in 
Propositions~\ref{prop-mon-ax-one} and \ref{prop-list}.  
The composite derived functors mentioned in
the introduction are then $\hh= U  L^{\mon} c R $ 
and $\pPsi = D c \phi^* N  Z c $ where $c$ denotes the 
cofibrant replacement functor in each category of monoids.  
Fibrant replacements are not necessary here because the right adjoints 
preserve all weak equivalences. 
\end{proof}

\begin{remark}
{\em An earlier version of this paper had a four step zig-zag of
Quillen equivalences instead of the above three step zig-zag. 
The earlier comparison of $\Spec(\ch)$ and $\Ch$ used 
two steps involving somewhat simpler and more naturally defined 
functors than $D$ and $R$.  
\[\xymatrix@=10mm{ \chZ \quad \ar@<.4ex>^-{F_0}[r] &
\quad  \spec(\chZ) \quad  \ar@<.4ex>^-{\text{Ev}_0}[l]
\ar@<-.4ex>_-{C_0}[r] &
\quad  \spec(\ch) \quad \ar@<-.4ex>_-{i}[l]
} \]
The inclusion of non-negatively graded chain complexes into $\mZ$-graded
chain complexes induces a functor $i\mc \spec(\ch) \to 
\spec(\chZ)$ with right adjoint $\co$, the prolongation of
the connective cover.  Then evaluation at the zeroth 
level gives a functor $\Ev\mc \spec(\chZ) \to \chZ$ with
left adjoint $F_0$~\cite[7.3]{hovey-stable}.
Here $\spec(\Ch)$ is $\Spec(\Ch, \bZ[1])$, the category of modules over 
Sym$(\bZ[1])$ in $\Ch^{\Sigma}$.

We show in Proposition~\ref{prop-extra} that
 $(\spec(\chZ), \spec(\ch), i, \co)$ and 
$ (\spec(\chZ),  \chZ, F_0, \Ev)$ are strong monoidal Quillen equivalences. 
The arguments for Proposition~\ref{prop-mon-ax-one} also extend to 
$\Spec(\Ch)$.  Thus, the arguments used in Theorem~\ref{thm-ch-hZ} would
also apply to give a four step zig-zag of equivalences involving
$\overline{\hh} = U L c \co f F_0 c$ 
and $\overline{\pPsi}  = \Ev f i \phi^* N Z c $ 
where $c$ and $f$ are the appropriate cofibrant and fibrant replacement
functors.  Here $U$, $\phi^*N$ and $i$ preserve all weak equivalences
so we have deleted the respective fibrant replacement functors. 
We should mention that with this longer zig-zag of functors 
we could only prove a non-natural version of Theorem~\ref{thm-Q-comm}.
}\end{remark}

\subsection{Extension to modules}
In this section we construct Quillen equivalences between the categories of 
modules over differential graded algebras and $HR$-algebra spectra.  Again,
we use the criteria developed in~\cite{emmc} to show that weak monoidal
Quillen equivalences induce equivalences on these categories of
modules.  The relevant statement,~\cite[3.12]{emmc}, requires one 
condition in addition to the Standing assumptions~\ref{stass} on the source 
and target model categories. 

\begin{definition}\label{def-Qinv}
{\em {\em Quillen invariance of modules} is said to hold in a monoidal model 
category if any weak equivalence of monoids $A \to B$ induces a Quillen 
equivalence between the associated model categories of modules
$A\Mod$ and $B\Mod$; see~\cite[3.11]{emmc} and~\cite[4.3]{ss}. 
}\end{definition}

As in the case of monoids, if 
$\lambda\mc \cD \adjoint \cC \mc R$ is a weak monoidal Quillen pair,
then for any monoid $A$ in $\cC$ the lax monoidal
right adjoint $R$ lifts to a functor 
$R \mc A\Mod \to RA\Mod$.   In~\cite[3.3]{emmc}, this module valued functor 
is shown to have a left adjoint $L^A\mc RA\Mod \to A\Mod$. 
Similarly, for $B$ a monoid in $\cD$, the
right adjoint $R$ induces a functor $(L^{\mon}B)\Mod \to B\Mod$
(which factors through $R(L^{\mon}B)\Mod$).   Again,~\cite[3.3]{emmc}
shows that there is a left adjoint denoted by $L_B\mc B\Mod \to 
(L^{\mon}B)\Mod$.  As with monoids, in general $L^A$ and $L_B$ do not
agree with $\lambda$ on underlying objects unless $(\lambda, R)$ is a strong
monoidal Quillen pair.  

\begin{theorem}\cite[3.12]{emmc}\label{thm-QI-mod}
Let $\lambda\mc \cD \adjoint \cC \mc R$ be a weak monoidal Quillen
equivalence such that $R$ preserves all weak equivalences.   Assume 
%Quillen invariance as well as 
the Standing assumptions~\ref{stass} 
hold for $\cC$ and $\cD$.  
\begin{enumerate} 
\item
For any cofibrant monoid $B$ in $\cD$ the adjoint pair
$L_B\mc B\Mod \adjoint (L^{\mon}B)\Mod \mc R$ is a Quillen equivalence.
\item
Suppose as well that Quillen invariance holds for $\cC$ and $\cD$. 
Then, for any monoid $A$ in $\cC$ the adjoint pair
$L^A\mc RA\Mod \adjoint A\Mod \mc R$ is a Quillen equivalence.
\end{enumerate} 
The analogues of these statements for modules over rings with many 
objects also hold; see~\cite[6.5]{emmc}.
\end{theorem}

The statement in~\cite[3.12]{emmc} is more general since it does
not require that the right adjoint $R$ preserves all weak equivalences.
Instead, one would require that $A$ is fibrant in (2) above.
In each of the monoidal Quillen equivalences considered in this
paper the right adjoint does preserve all weak equivalences though.
We should remark that here (2) follows from (1) for $B= c(RA)$ by 
using Quillen invariance since $A\to L^{\mon}c(RA)$ is a weak
equivalence by Theorem~\ref{thm-emmc}. 

The following proposition is proved as Corollary~\ref{cor-mon-ax}. 

\begin{proposition}\label{prop-QI}
Quillen invariance for modules holds in $H\bZ\Mod$, 
$\spec(\sab)$, $\spec(\ch)$ and $\chZ$.
\end{proposition}

Equivalences of the associated categories of modules then follow from
Theorem~\ref{thm-QI-mod}. 

\begin{corollary}\label{Big Cor}
Using the composite functors $\hh\mc \dga_{\mZ} \to H\mZ\alg$ 
and $\pPsi \mc H\mZ\alg \to \dga_{\mZ}$ 
defined above we have the following equivalences of module categories. 
\begin{enumerate}
\item\label{dos} For an (associative) differential graded algebra $A$,
there is a Quillen equivalence between differential graded $A$-modules and
$\hh A$-module spectra.
\[A\Mod \simeq_Q \hh A\Mod\] 

\item\label{tres} For $B$ an (associative) $H\mZ$-algebra spectrum, there is a 
Quillen equivalence between differential graded $\pPsi B$-modules and
$B$-module spectra.
\[ \pPsi B\Mod \simeq_Q B\Mod\]
\end{enumerate}
The analogues of these statements for modules over rings with many
objects also hold.
\end{corollary}

\begin{proof}
%The functor $\Phi$ in Theorem~\ref{thm-ch-hZ}, Part~\ref{dos} 
%is the composite functor
%$\Phi= U f L c R f$ where $c$ and $f$ are the appropriate
%cofibrant and fibrant replacement functors of monoids.  
%These cofibrant and fibrant replacement functors induce
%Quillen equivalences by the Quillen invariance property for modules over
%weakly equivalent monoids~\cite[4.3]{ss} which is one of the hypotheses
%on the categories involved in an applicable Quillen equivalence.
Both parts follow from repeated applications of Theorem~\ref{thm-QI-mod}
and Quillen invariance statements.  
%
%The functor $\pPsi$ in Theorem~\ref{thm-ch-hZ}, Part~\ref{tres} is 
%the composite functor
%$\pPsi = D c \phi^* N f Z c $ with $c$ and $f$ as above.
%Again Part~\ref{tres} follows from repeated applications
%of Theorem 3.12, Parts 1 and 2 from~\cite{emmc}
%and Quillen invariance statements.   Since $\phi^* N$
%preserves all weak equivalences the fibrant replacement could be deleted. 
%In the proof of Theorem~\ref{thm-Q-comm} it is particularly important that
%we have this simplified version, $\pPsi' = \Ev f i \phi^* N Z c$.
%
The statements 
%similar to Parts~\ref{dos} and~\ref{tres} hold 
for modules over enriched categories (or rings with many objects) follow 
from Theorem 6.5, Parts 1 and 2 from~\cite{emmc}.  
\end{proof}

We use the Quillen equivalences of modules over $H\bZ$-algebras and
differential graded algebras from Corollary~\ref{Big Cor} to construct 
algebraic models for
{\em rational stable model categories},  
the  pointed model categories $\C$ where suspension is an invertible functor 
on the homotopy category and for any two objects $X, Y$ the homotopy
classes of maps $[X,Y]^{\Ho\C}$ form a rational vector space. 
We show that basically any rational stable model category with 
a %set of 
small generator %s 
is Quillen equivalent to (right) 
modules over a rational differential graded algebra. 
%with many objects (a rational DG-category).
This is the rational version of~\cite[3.1.1]{stable} which shows that basically
any stable model category $\C$ with a %set of 
small generator %s 
$\G$ is Quillen equivalent
to (right) modules over a ring spectrum %with many objects, 
$\E(\G)$, the endomorphism ring of $\G$. 
We must say ``basically" here because in fact we need $\C$ to
be Quillen equivalent to a spectral model category.  A {\em spectral
model category} is the analogue of a simplicial model category
where simplicial sets have been replaced by symmetric spectra (or
one of the other monoidal model categories of spectra); 
see~\cite[3.5.1]{stable} and Proposition~\ref{prop-spectral} below. 
We now state the more general case where we allow a set of compact
generators, this is based on~\cite[3.3.3, 3.9.3]{stable}.

\begin{corollary}\label{cor-rational}
Let $\C$ be a 
%cofibrantly generated, proper, simplicial, rational 
rational stable model category with a set $\G$ of small generators which
is Quillen equivalent to a spectral model category.  Then
there exists a differential graded algebra with many objects $\cA$
%\iso \pPsi (H\mQ \sm c\E(\G))$ 
and a chain of
Quillen equivalences between $\C$ and the category of right $\A$-modules.
The objects of $\A$ correspond to the objects in $\G$ and 
there is an isomorphism of graded $\mQ$-categories
between the homology category $H_*\cA$ and
the full graded subcategory of $\Ho(\C)$ with objects $\G$.
\[ \C \simeq_Q \Modr \E(\G) \simeq_Q \Modr \A \]
\end{corollary}
\noindent
To specify $\cA$, recall the functor $\pPsi$ from 
Theorem~\ref{thm-ch-hZ}.  This functor can be
extended to rings with many objects; see~\cite[6]{emmc} 
and~\cite[6]{add}.
Then $\cA \iso \pPsi (H\mQ \sm c\E(\G))$ where $c$ is a cofibrant
replacement functor for ring spectra (with many objects) from~\cite[6.3]{emmc}.

In general, $\Modr\A$ is not a practical algebraic model; it is difficult to 
make explicit and $\cA$ is quite large.  This model can be used as
a stepping stone to a practical model though.
For example, this corollary applies to the category of rational $G$-equivariant
spectra for any compact Lie group $G$.  
In~\cite{shi-circle} and~\cite{gs} this large model is used
to show that there is an explicit algebraic model 
%given in~\cite{g-circle} and~\cite{g-tori} are 
Quillen equivalent to the category of rational $T^n$-equivariant spectra for 
$T^n$ the $n$-dimensional
torus.  For $G$ finite this also extends the results of~\cite[Appendix A]{gm}
to incomplete universes; see~\cite[5.1.2]{stable}.
On the other hand, although this corollary applies to the rational motivic 
stable homotopy theory of schemes from~\cite{jardine, voevodsky} 
no such explicit algebraic model is known.    

Several different sets of hypotheses can be used to 
ensure that $\C$ is Quillen equivalent to a spectral model category. 
A partial converse from~\cite[3.5.2]{stable} shows that any spectral
model category is stable and simplicial.

\begin{proposition}\label{prop-spectral}
If one of the following lists of hypotheses holds, then $\C$ is
Quillen equivalent to a spectral model category.
\begin{enumerate}
\item \cite[1.8, 4.7a, 4.9a]{dugger-spectral} $\C$ is stable and combinatorial.
\item \cite[8.11, 9.1]{hovey-stable} $\C$ is stable, simplicial, left proper, cellular and the domains of the generating cofibrations are cofibrant.
\item \cite[3.8.2]{stable} $\C$ is stable, simplicial, proper and cofibrantly generated.
\end{enumerate}
\end{proposition}

\begin{proof}[Proof of Corollary~\ref{cor-rational}]
If $\C$ is a rational stable model category with a set of small 
generators which is Quillen equivalent to a spectral model category $\D$,
then it follows that $\D$ is also rational and has a set $\G$ of small 
generators.  In~\cite[3.9.3]{stable} we show that any spectral model 
category $\D$ with a set $\G$ of small generators is equivalent
to (right) modules over $\E(\G)$, the endomorphism ring spectrum of 
$\G$ (a symmetric ring spectrum with many objects) defined 
in~\cite[3.7.5]{stable}.  Since there is a natural
isomorphism of the homotopy groups of $\E(\G)$ and the graded homotopy
classes of maps between elements in $\G$ by~\cite[3.5.2]{stable}, 
the homotopy groups of $\E(\G)$ are rational when $\D$ is rational. 

If $R$ is a ring spectrum with rational homotopy groups which is 
cofibrant as an underlying module, then the
map $R \to H\mQ \sm R$ is a stable equivalence of ring spectra.  
Since $\mQ$ is flat over $\pi_*^s$, this
follows from the spectral sequence for computing the homotopy 
of $H\mQ \sm R$~\cite[IV.4.1]{ekmm} and the Quillen equivalences of
ring spectra established in~\cite{sch}. 
%and~\cite[5.3.8, 5.3.10]{hss} 
%since cofibrant symmetric ring spectra are also 
%cofibrant as underlying spectra~\cite[?]{ss}.  
Similarly, if $\cR$ is a ring spectrum with many objects with
rational homotopy groups such that each $\cR(i,j)$ is cofibrant
as a module, then $\cR \to H\bQ \sm \cR$ is a stable equivalence
of ring spectra with many objects.  

Quillen invariance then shows that $\Modr \E(\G)$, $\Modr c\E(\G)$ and 
$\Modr (H\mQ \sm c\E(\G))$  are Quillen equivalent, where $c$ denotes the 
cofibrant replacement of ring spectra (with many objects) 
from~\cite[6.3]{emmc}.  Here we are using the fact that, since the sphere
spectrum is cofibrant, the cofibrant replacement is pointwise cofibrant
as a module. Then $\Modr (H\mQ \sm c\E(\G))$ is Quillen
equivalent to the model category of differential graded 
$\cA = \pPsi (H\mQ \sm 
c\E(\G))$-modules, $\Modr \cA$, by the many objects analogue
of Corollary~\ref{Big Cor}, Part~\ref{tres}.  This analogue
holds by Theorem 6.5, Parts 1 and 2 
from~\cite{emmc}
applied to each of the three adjoint pairs in Proposition~\ref{prop-list}.
%The necessary cofibrant replacements for ring spectra
%with many objects follow from Theorem 6.3 from~\cite{emmc}. 

%Since $\hh \mZ [0]$ is stably equivalent to $H\mZ$ and 
Since $\pPsi$ and $\hh$ 
induce an equivalence between the homotopy
categories of differential graded algebras and $H\mZ$-algebras, 
the calculation of $H_* \cA = H_* \pPsi (H\mQ \sm c\E(\G))$ follows 
from~\cite[3.5.2]{stable},
the stable equivalence of $\E(\G)$ and $H\mQ\sm c \E(\G)$ and
the stable equivalence of $\hh \mZ [0]$ and $H\mZ$.   This last fact
follows since $\mZ [0]$ is the unit of the monoidal structure on $\chZ$ and 
each of the functors $U$, $L^{\mon}$ and $R$ preserve these units.  This is
simple to verify for $U$ and $R$ and is verified for $L^{\mon}$ in the proof of
Proposition~\ref{prop-sab-ch}; see also~\cite[3.7]{emmc}.  
\end{proof}

%%%%%%%%%%%%%%%%%  I think all of this is true, but so much technicalities..
%%%%% need to show \tQ preserves all weak equivalences, and \tQ left exact...
%%%%% other way is easier.............
%The functor $\tQ\mc \cS_* \to s\V$ with $\tQ K$ the reduced free simplicial 
%rational vector space on the pointed simplicial set $K$ is symmetric monoidal.
%Thus applying $\tQ$ levelwise takes symmetric ring spectra (with
%many objects) to $H\mQ$-algebras (enriched categories over $H\mQ$-modules) 
%since $\tQ(F_0S^0) \iso H\mQ$.  Moreover, for symmetric ring
%spectra (with many objects) with rational homotopy groups the 
%inclusion $i \mc R \to \tQ R$ is a stable equivalence.   To see this, first
%note that the map $H\mQ \to \tQ H\mQ$ induces an isomorphism of stable
%homotopy groups.  Since $\tQ$ is left exact and $H\mQ$ generates
%all $H\mQ$-modules, it follows that $i$ is a stable equivalence for any 
%cofibrant $H\mQ$-module.     
%But any $H\mQ$-module is level equivalent
%to a cofibrant one and $\tQ$ preserves level equivalences, so $i$ is
%a stable equivalence for any $H\mQ$-module.  Finally, for any spectrum $X$
%with rational homotopy groups the  
% $H\mQ$ generates the category of rational spectra and
%$\tQ$ is left exact.  Notice, $H\mQ$ obviously generates the category
%of $H\mQ$-modules, but any rational spectrum is weakly equivalent to
%an $H\mQ$-module since the map $X \to H\mQ \sm X$ .  

\section{Monoidal model categories and the monoid axiom}\label{sec-monoid}

In this section we prove Proposition~\ref{prop-mon-ax-one}
which sates that the Standing assumptions~\ref{stass} hold for each of the
four categories $H\bZ\Mod$, $\spec(\sab)$, $\spec(\ch)$ and $\Ch$.  
Recall, the standing assumptions require that these categories are
cofibrantly generated, monoidal model categories (see Definition~\ref{def-mmc})
which satisfy the monoid axiom (see Definition~\ref{def-ma}), have cofibrant 
unit objects and have all objects small. 
We also establish Proposition~\ref{prop-QI} by showing that Quillen 
invariance (see Definition~\ref{def-Qinv}) holds in each of these four 
categories.  

The proof that all objects are small in each of these four categories
is delayed to the end of this section in Propositions~\ref{prop-loc}
and~\ref{prop-small}.  All of the other statements for $H\bZ$-module 
spectra follow from~\cite[3.4.9, 5.3.8, 5.4.2, 5.4.4]{hss} and the
fact that $H\bZ$ is cofibrant in $H\bZ\Mod$.  
We are left with establishing the standing assumptions and Quillen invariance 
for the other three categories. 
%$\spec(\sab)$ $\spec(\ch)$. 
%see Definition~\ref{def-Qinv}.  

We first verify that $\spec(\sab)$, $\spec(\ch)$ and $\Ch$ are 
cofibrantly generated, monoidal model categories.
$\chZ$ is shown to be a cofibrantly generated, monoidal model category 
in~\cite[2.3.11, 4.2.13]{hovey-book} and the proof obviously applies to 
$\ch$ as well.  These properties follow for $\sab$ from~\cite[SM7]{Q} 
by using the strong symmetric monoidal functor $\tZ$ from pointed simplicial 
sets to $\sab$; see also~\cite{emmc}. 
Since $\sab$ and $\ch$ are cofibrantly generated, monoidal
model categories, the same properties for $\spec(\sab)$ and 
$\spec(\ch)$ then follow from~\cite[8.11]{hovey-stable}.   Here we
are also using the fact that  the generating cofibrations in both 
$\sab$ and $\ch$ are maps between cofibrant objects.  

\cite{hovey-stable} was not able to verify in general though
that $\spec(\C)$ satisfies the monoid axiom if $\C$ does.
Instead we follow an approach for verifying the monoid axiom 
for $\spec(\sab)$, $\spec(\ch)$ and $\Ch$  
modeled on the approach for $\spec(\sset)$ in~\cite[5.4]{hss}.  To establish 
the general pattern we consider $\Ch$ first.  

We prove the monoid axiom for $\Ch$ by using both the projective and
injective model structures.  The {\em projective model structure} is the
structure we have referred to so far in this paper; 
see~\cite[2.3.11]{hovey-book}.  The weak equivalences are the 
quasi-isomorphisms, the (projective) fibrations are the surjections and the 
(projective) cofibrations 
are the maps with the left lifting property with respect to the trivial 
fibrations.   We now introduce the {\em injective model structure}
from~\cite[2.3.13]{hovey-book}.  The weak equivalences  are again the
quasi-isomorphisms, the injective cofibrations are the monomorphisms
and the injective fibrations are the maps with the right lifting
property with respect to the trivial cofibrations. 
 
\begin{proposition}\label{prop-ma-ch} 
The monoid axiom holds for $\Ch$.
\end{proposition} 

\begin{proof}
We first consider the monoid axiom for the generating trivial 
(projective) cofibrations.  Namely, we show that 
any map obtained by pushouts and compositions of maps of the
form $j \otimes \id_Z\mc A \otimes Z \to B \otimes Z$ is a weak
equivalence where  
$j: A \to B$ is any generating trivial (projective) cofibration and $Z$ is any 
object.  It then follows by~\cite[3.5(2)]{ss} 
that the full monoid axiom holds.   The generating trivial 
(projective) cofibrations are all of the form 
$j\mc 0 \to D^n$~\cite[2.3.3]{hovey-book}. 
Here $D^n$ is the acyclic chain complex with all degrees zero
except $(D^n)_n = \bZ$ and $(D^n)_{n-1} = \bZ$.  One can show that 
$D^n \otimes Z$ is acyclic for any $Z$ and thus
$j \otimes \id_Z \mc 0 \to D^n \otimes Z$ is a monomorphism 
and a quasi-isomorphism.  The monoid axiom then follows since
maps which are monomorphisms and quasi-isomorphisms are preserved under 
pushouts and compositions.
This last statement follows, for example, from the fact that the trivial 
cofibrations in the injective model structure on $\Ch$ are the 
maps which are monomorphisms and quasi-isomorphisms. 
\end{proof}

As with $\Ch$, we establish the monoid axiom on $\spec(\sab)$ and $\spec(\ch)$ 
by forming {\em injective stable model structures.}  We then refer to the 
stable model structures introduced above from~\cite{hovey-stable} as
the projective stable model structure.  Throughout this paper we use
{\em `injective'} to refer to this new model structure;  so 
`injective' implies certain lifting properties, and not necessarily 
`monomorphism'.   

\begin{proposition}\label{prop-inj-MC}
There is an injective stable model structure on both of the categories 
$\spec(\sab)$ and $ \spec(\ch)$ 
%and $\spec(\chZ)$ 
with injective cofibrations the monomorphisms,
weak equivalences the (projective) stable equivalences from
\cite{hovey-stable} and fibrations the maps with the right lifting
property with respect to the trivial injective cofibrations.
Note, here all objects are injective cofibrant.
\end{proposition}

The reason these injective stable model categories are so useful is
that there is an action of the respective projective stable model categories
on these injective model categories which satisfies an analogue of the
pushout product axiom.  This action appeared implicitly in the proof of
Proposition~\ref{prop-ma-ch} to show that the tensor of a trivial projective
cofibration and an injective cofibrant object is a trivial injective 
cofibration.  This action is analogous to the structure of simplicial 
model categories which have an action by simplicial sets.  
As with the definition of monoidal model categories, the unital
condition is automatic and thus suppressed because the units here
are all cofibrant.

\begin{proposition}\label{prop-inj-proj}
If $i\mc A \to B$ is an injective cofibration in 
$\spec(\sab)$, $\spec(\ch)$ or $\Ch$ 
% and $\spec(\chZ)$ 
and $i'\mc K \to L$ is a projective cofibration 
on the same underlying category then 
   $$i\boxprod i'\mc A \otimes L   \coprod_{  A\otimes K}    B\otimes   K \to B \otimes L$$
is an injective cofibration which is a weak equivalence if either one
of $i$ or $i'$ is.
That is, the injective stable model structure in each case
is a {\em Quillen module} over the projective stable model structure; see
also~\cite[4.2.18]{hovey-book}.
\end{proposition}

We prove this proposition for $\Ch$ first and
delay the more involved proofs of 
Propositions~\ref{prop-inj-MC} and~\ref{prop-inj-proj} for
$\spec(\sab)$ and $\spec(\ch)$ until after Corollary~\ref{cor-mon-ax}
below.

\begin{proof}[Proof of Proposition~\ref{prop-inj-proj} for $\Ch$]
Recall the injective and projective model structures on $\Ch$
defined above Proposition~\ref{prop-ma-ch}.
%This is not much of a reduction for the injective cofibrations since the set 
%$I'$ of generating injective cofibrations
%from~\cite[2.3.14] contains a map from each isomorphism class of monomorphisms
%between countable chain complexes.   (Here the cardinality of a chain
%complex is given by the cardinality of the union of all the levels.) 
We use the criterion in~\cite[3.5]{ss} to reduce to checking the 
pushout product axiom for $K\to L$ a generating (trivial) projective cofibration
and $A \to B$ any (trivial) monomorphism (that is, a (trivial) injective
cofibration). 

The generating set of projective cofibrations is given by 
$I = \{i_n\mc \bZ[n-1] \to D^n\}$.
Here $\bZ[n]$ is the complex concentrated in degree $n$ 
%with $(S^n)_n = \bZ$
and $D^n$ is the acyclic complex modeling the $n$-disk. 
As graded abelian groups (ignoring differentials) the pushout product of 
$i\mc A \to B$ and $i_n$ is the map $i \boxprod i_n \mc
(B \otimes \bZ[n-1]) \oplus (A \otimes \bZ[n]) 
\to (B \otimes \bZ[n-1]) \oplus (B \otimes \bZ[n])$.  This is a (trivial) injective
cofibration  if $i\mc A \to B$ is.   Also, for 
$j_n \mc 0 \to D^n$, $i \boxprod j_n \mc A \otimes D^n \to B \otimes D^n$
is a map between acyclic complexes and hence is a quasi-isomorphism.       
%If in addition one of these maps is a quasi-isomorphism, then
%so is the pushout product.  This follows from applying the long exact sequence 
%in homology to maps of cofiber sequences and the fact that if $P$ is
%projective cofibrant then $P \otimes A$ is acyclic if either $P$ or $A$ is. 
\end{proof}

Next we deduce the
monoid axiom and the Quillen invariance properties for
$\spec(\sab)$, $\spec(\ch)$ and $\Ch$ assuming the Quillen module
structure from Proposition~\ref{prop-inj-proj}. 
Propositions~\ref{prop-mon-ax-one} and~\ref{prop-QI} then follow for these 
three categories since the unit objects are all cofibrant and all objects are 
shown to be small in Propositions~\ref{prop-loc} and~\ref{prop-small}.
As mentioned above, the statements about $H\bZ$-modules follow 
from~\cite[5.4.2, 5.4.4]{hss}.  

\begin{corollary}\label{cor-mon-ax}
The stable model structures on $\spec(\sab)$, $\spec(\ch)$ and $\Ch$
are cofibrantly generated, monoidal model 
categories which satisfy the monoid axiom.  
Quillen invariance for modules also holds in each of these 
categories.  
\end{corollary}

\begin{proof}
We verified that these model structures are cofibrantly generated,
monoidal model categories at the beginning of this section.
We use Proposition~\ref{prop-inj-proj} to establish
the monoid axiom.  Since the injective cofibrations are monomorphisms, any 
object $A$ is injective cofibrant. 
Hence, $A \sm -$ takes stably trivial projective cofibrations to 
stably trivial injective cofibrations by Proposition~\ref{prop-inj-proj}.  
Pushouts and directed colimits of 
stably trivial injective cofibrations are again stably trivial injective
cofibrations.
Since the injective stable equivalences agree with the projective stable
equivalences this establishes the monoid axiom.

In~\cite[6.2(c)]{add}, Quillen invariance 
is shown to follow from the following two properties (see also~\cite[4.3]{ss}):
\begin{itemize}
\item[(QI1)]
For any projective cofibrant object $A$ and any weak 
equivalence $X\ra Y$, the map $A\tens X \ra A\tens Y$ is also a weak
equivalence.
\item[(QI2)]
Suppose $A\cofib B$ is a projective cofibration, and $X$ is any object.  
Then for any map $A\tens X \ra Z$, the map from the homotopy pushout of 
$B\tens X \lla A\tens X \ra Z$ to the pushout is a weak equivalence.
\end{itemize}

These two properties in turn follow from Proposition~\ref{prop-inj-proj}.  
For $A$ a projective
cofibrant object, $A \sm -$ preserves stably trivial injective 
cofibrations by Proposition~\ref{prop-inj-proj}.  Hence, by Ken Brown's 
lemma~\cite[1.1.12]{hovey-book}, 
$A \sm -$ preserves stable equivalences between injective cofibrant objects.
Since all objects are injective cofibrant, this implies (QI1). 

By Proposition~\ref{prop-inj-proj}, $A \tens X \to B \tens X$ is an injective 
cofibration whenever $A \cofib B$ is a projective cofibration
since every object $X$ is injective cofibrant.   (QI2) then follows from the 
fact that the injective stable model structure is left proper
since all objects are cofibrant; see also~\cite[5.2.6]{hovey-book}.
\end{proof}

We now turn to the proofs of Propositions~\ref{prop-inj-MC} 
and~\ref{prop-inj-proj}.   
The proofs of these propositions are very similar for $\spec(\sab)$ and 
$\spec(\ch)$.  We consider in detail the case of $\spec(\ch)$.  
For $\spec(\sab)$, one can follow the same outline at every stage;
see Proposition~\ref{prop-cof-gen} for an identification of
the generating trivial cofibrations.

We will construct the injective stable model categories as localizations
of injective level model categories.   This follows the outline
in~\cite{hovey-stable} which first constructs 
projective level model categories  
on $\spec(\sab)$ and $\spec(\ch)$ with weak equivalences the level
equivalences, fibrations the level fibrations and cofibrations the 
maps with the left lifting property with respect to the level trivial
fibrations.  These are cofibrantly generated, monoidal model categories 
by~\cite[8.3]{hovey-stable}.  The next lemma establishes the 
analogous injective level model categories.

\begin{lemma}\label{prop-inj-level}
There is an injective level model category on both of the categories 
$\spec(\sab)$ and $\spec(\ch)$ 
with cofibrations the monomorphisms, weak equivalences the level equivalences 
and fibrations the maps with the right lifting
property with respect to the level trivial cofibrations.
\end{lemma}

\begin{proof}
First consider the injective model structure
on $\ch$ with cofibrations the monomorphisms and weak equivalences the
quasi-isomorphisms;~\cite[2.3.13]{hovey-book} discusses the analogue on $\chZ$.
This model category is cofibrantly generated; a set of generating 
cofibrations $I'$ contains a map from each isomorphism class of monomorphisms 
between countable chain complexes and $J'$ is the set of weak equivalences
in $I'$.  Here the cardinality of a chain
complex is given by the cardinality of the union of all of the levels.   

The proof that the injective level model structure on $\spec(\ch)$ 
forms a cofibrantly generated model category follows as in~\cite[5.1.2]{hss}.
The only changes needed are that one uses a chain homotopy instead of
a simplicial homotopy in~\cite[5.1.4(5)]{hss} and for~\cite[5.1.6, 5.1.7]{hss} 
the arguments for simplicial sets are replaced by the analogues 
of~\cite[2.3.15, 2.3.21]{hovey-book} for $\ch$ instead of $\chZ$ as above.
(For $\spec(\sab)$ the proof follows~\cite[5.1.2]{hss} almost verbatim.) 
Note that every spectrum is cofibrant here.  Again 
the generating cofibrations are the monomorphisms 
between countable objects and the set of generating trivial cofibrations is  
the subset of level weak equivalences.  Call the fibrant objects
here the {\em injective fibrant objects}.

One could also establish the injective level model structures by
using Jeff Smith's work on combinatorial model categories. %see~\cite{bekeI}.  
We sketch the verification of the criteria listed in~\cite[1.7]{bekeI}. 
Criteria (c0) and (c3) follow from~\cite[1.15, 1.18, 1.19]{bekeI} since
the functors $ev_n\mc \spec(\ch) \to \ch$ are accessible and
$\ch$ is locally presentable by Proposition~\ref{prop-loc}; see 
also~\cite[7.3]{dugger-pres} and \cite[2.23, 2.37]{adamek}.
Maps with the right lifting property with respect to monomorphisms must be  
levelwise weak equivalences
(and levelwise fibrations) since the generating projective cofibrations are
monomorphisms~\cite[8.2]{hovey-stable}; criterion (c1) thus holds.
Finally, directed colimits and pushouts preserve levelwise weak equivalences
which are monomorphisms since colimits are created on each level; this
verifies criterion (c2). 
\end{proof}

We next show that the analogue of Proposition~\ref{prop-inj-proj} holds 
for the associated level model categories. 

\begin{lemma}\label{lem-proj-level}
The injective and projective level model categories on 
$\spec(\sab)$ and  $\spec(\ch)$ 
satisfy the property stated in Proposition~\ref{prop-inj-proj}. 
That is, the injective level model structure in each case
is a Quillen module over the projective level model 
structure~\cite[4.2.18]{hovey-book}.
\end{lemma}

\begin{proof}
Again, use~\cite[3.5]{ss} to reduce to considering only the generating
(trivial) projective cofibrations.  Each part of the proof that the 
injective level model structure is a Quillen module over the projective level 
model structure then follows from checking that the injective model structure 
on $\ch$ is a Quillen module over the projective
model structure on $\ch$~\cite[4.2.18]{hovey-book}.  This follows for $\ch$ as
outlined in the proof of Proposition~\ref{prop-inj-proj} for $\Ch$ above.
\end{proof}

\begin{proof}[Proof of Proposition~\ref{prop-inj-MC}]
The injective level model structure is left proper and cellular.  Thus, 
we may localize with respect to the set of maps $\J= \{F_{n+1}(A \otimes \mZ[1])
\to F_nA\}_{n \geq 0}$ for one object $A$ from each isomorphism class of countable objects 
in $\spec(\ch)$.  See~\cite[Section 2]{hovey-stable} for a brief summary
of Bousfield localization; the definitive reference is~\cite{psh}.
One could also use the machinery developed by Jeff Smith~\cite{smith} to 
establish
the existence of these local model categories because these categories
are also left proper and combinatorial (cofibrantly generated and locally
presentable, see Proposition~\ref{prop-loc}); see also~\cite[Section 2]{dugger}.

We are left with showing that the 
injective stable (local) equivalences agree with the (projective) stable 
equivalences defined in~\cite[8.7]{hovey-stable} for the projective stable
model structure on $\ch$.  
Following~\cite[8.8]{hovey-stable} one can
show that the $\J$-local objects, or injective stably fibrant objects are the 
injective fibrant objects which are also 
$\Omega$-spectra~\cite[8.6]{hovey-stable}.  
Since injective fibrations are the maps
with the right lifting property with respect to all level trivial
monomorphisms, they are in particular level fibrations.    
Any injective fibrant $\Omega$-spectrum is thus a levelwise fibrant
$\Omega$-spectrum, that is, a projective fibrant $\Omega$-spectrum.  

By definition, a map $f$ is an injective stable (local) equivalence if 
and only if $\map_{\inj}(f, E)$ is a weak equivalence 
for any injective fibrant $\Omega$-spectrum $E$.
Here $\map_{\inj}(f,E)$ is the homotopy function complex in the injective level 
model structure as constructed in~\cite[Section 2]{hovey-stable}, 
\cite[5.4]{hovey-book} or~\cite[17.4.1]{psh}.  
Similarly, a map $f$ is a (projective) stable equivalence if and only
if $\map_{\proj}(Qf, F)$ is a weak equivalence for any projective fibrant
$\Omega$-spectrum $F$.  Here $Q$ is the projective cofibrant
replacement functor.  

The identity functors induce a Quillen equivalence between the injective level
model structure and the projective level model structure since weak 
equivalences in both cases are level equivalences.  Thus, 
by~\cite[17.4.16]{psh}, if $X$ is projective cofibrant and $Y$ is
injective fibrant then $\map_{\inj}(X, Y)$ is weakly equivalent to
$\map_{\proj}(X,Y)$. 
As above, denote by $Q$ the projective cofibrant replacement 
functor.   Since $QX \to X$ is a level equivalence, $f$ is an injective 
stable equivalence if and only if $\map_{\inj}(Qf, E)$,
or equivalently $\map_{\proj}(Qf, E)$,
is a weak equivalence for any injective fibrant $\Omega$-spectrum $E$.
It follows that if $f$ is a projective stable equivalence, then $f$ is
also an injective stable equivalence since injective fibrant 
$\Omega$-spectra are also projective fibrant $\Omega$-spectra.
Similarly, denote by $R$ the injective fibrant replacement functor.
Since $Y \to RY$ is a level equivalence, $f$ is a projective
stable equivalence if and only if $\map_{\proj}(Qf, RF)$,
or equivalently $\map_{\inj}(Qf, RF)$,
is a weak equivalence for any projective fibrant $\Omega$-spectrum $F$.
It follows that if $f$ is an injective stable equivalence, then $f$ is
also a projective stable equivalence. 
\end{proof}

\begin{proof}[Proof of Proposition~\ref{prop-inj-proj}]
We now show that the injective stable model structure is a
Quillen module over the projective stable model structure.  The pushout
product of an injective cofibration and a projective cofibration
is an injective cofibration by Lemma~\ref{lem-proj-level}. 

We next 
consider the case where $j\mc K \to L$ is a stably trivial injective
cofibration and $i\mc F_k(B) \to F_k(C)$ is a generating projective cofibration.
The functor $- \otimes_{\Sym(\mZ[1])} F_k(B) = - \sm F_k B$ preserves
injective cofibrations by Lemma~\ref{lem-proj-level}.  It also
preserves stably trivial injective cofibrations for $B$ a countable projective 
cofibrant complex in $\ch$.  This follows by the construction of
the injective stable model category since 
this functor takes maps in $\J$ into $\J$ up to isomorphism; 
see the last sentence of~\cite[2.2]{hovey-stable}.  
In fact, for
each map in $\J$ it just replaces $A$ by $A \otimes B$ and $n$ by $n +k$.  
Since $B$ and $C$ are countable complexes, $j \sm F_k(B)$ and   
$j \sm F_k(C)$ are both stably trivial injective cofibrations.  
Let $P = (L \sm F_k(B))\sm_{K \sm F_k (B)} K \sm F_k(C)$; 
it follows that $j' \mc K \sm F_k(B) \to P$ is a stably trivial injective 
cofibration since it is the pushout of $j \sm F_k(A)$.  
The composite of $j'$
and $j \boxprod i \mc P \to L \sm F_k(B)$ is $j \sm F_k(B)$, so
by the two out of three property, $j \boxprod i$ is a stable 
equivalence as well.   
By the criterion in~\cite[3.5]{ss}, it follows that this holds for  
$i$ any projective cofibration. 

We next consider the case where $j\mc K \to L$ is an injective 
cofibration and $i\mc A \to B$ is a stably trivial projective cofibration.
In the next paragraph we show that $K\sm -$ takes stably trivial projective 
cofibrations to stably trivial injective cofibrations. 
It then follows by the same argument as in the last paragraph that the 
pushout product of an 
injective cofibration and a stably trivial projective cofibration is a
stably trivial injective cofibration.  

Since $K$ is injective cofibrant for any object, $K \sm -$ takes (level trivial)
projective cofibrations to (level trivial) injective cofibrations by 
Lemma~\ref{lem-proj-level}.  
So to show $K \sm -$ takes stably trivial projective cofibrations
to stably trivial injective cofibrations, it is enough to show that
$K \sm -$ takes stably trivial projective cofibrations to stable equivalences. 
By~\cite[3.5]{ss} it is enough to consider $K \sm -$ on the 
generating stably trivial projective cofibrations. 
Since the generating stably trivial projective cofibrations
are maps between projective cofibrant objects, we 
reduce our problem to showing that $K \sm Q(-)$ preserves stable 
equivalences where $Q$ is a projective cofibrant replacement functor.   

From Lemma~\ref{lem-proj-level} we see that 
for a projective cofibrant object $QX$, the functor $ - \otimes QX$
preserves level trivial injective 
cofibrations.  Thus, by Ken Brown's lemma~\cite[1.1.12]{hovey-book}, 
$- \otimes QX$ also preserves level equivalences between injective cofibrant
objects, i.e. all level equivalences.  Since $QK \to K$ is
a level equivalence, it follows that $K \otimes Q(-)$ is level equivalent to 
$QK \otimes Q(-)$. So we only need to see that $QK \sm Q(-)$ preserves
stable equivalences. 
Since the projective stable model category is a monoidal model category
by~\cite[8.11]{hovey-stable}, $QK \sm -$  preserves projective stably
trivial cofibrations.  Then, by Ken Brown's lemma again,
$QK \sm -$ preserves stable equivalences between projective cofibrant
objects;  that is, $QK \sm Q(-)$ preserves stable equivalences. 
\end{proof}

Finally, we verify the requirement that each object is small
by showing that these categories are locally presentable.  
%Any object in a locally presentable category is small;  see  
%~\cite[1.13, 1.16, 1.17]{adamek} or~\cite[5.2.10]{bor}.
The small object argument, see~\cite[2.1.14]{hovey-book} or~\cite[7.12]{ds},  
then applies to any set of maps
and is used in~\cite{ss}
to construct the model structures on the categories of algebras and modules. 

\begin{proposition}\label{prop-loc}
Each of the categories 
$H\mZ\Mod, \spec(\sab), \spec(\ch)$  and $\chZ$  
is locally presentable. 
The categories of monoids, modules and algebras over these
categories are also locally presentable. 
\end{proposition}

\begin{proof}
Recall $D^n$ is the acyclic chain complex modeling the $n$-disk.
$\chZ$ is locally presentable by~\cite[1.20]{adamek} because
the set $\{ D^n \}_{n \in \mZ}$ is a strong generator with each
object finitely presentable; see~\cite[0.6, 1.1]{adamek}.  
Similar arguments apply to $\ch$.
Also, the category of simplicial
sets, $\simpu$, is locally presentable by~\cite[5.2.2b]{bor} because it is 
set-valued diagrams over a small category.  

~\cite[5.5.9]{bor} and~\cite[5.3.3, 5.7.5]{bor} show that 
over a locally presentable category any category of algebras over a monad 
which commutes with directed colimits or any functor category
from a small category to a locally presentable category is again
a locally presentable category.  Since  
$H\mZ\Mod, \spec(\sab)$ and $\spec(\ch)$   
can be built using these two methods from $\ch$ or $\simpu$,
the statement follows.
\end{proof}

\begin{proposition}\label{prop-small}
\cite[1.13, 1.16, 1.17]{adamek},~\cite[5.2.10]{bor}
Each object in a locally presentable category is small relative to
the whole category.
\end{proposition}

\section{Proof of Proposition~\ref{prop-list}} \label{sec-proof}

Before considering each of the three parts of Proposition~\ref{prop-list}
separately, we identify the generating cofibrations and
trivial cofibrations in $\spec(\sab)$.  This is useful in the two parts
involving this model category.  
Recall that $\tZ\mc \simp \to \sab$ is the reduced free abelian group functor.

\begin{proposition}\label{prop-cof-gen}
The generating cofibrations and trivial cofibrations in the 
cofibrantly generated stable model structure on $\spec(\sab)$
defined in~\cite{hovey-stable} are given by applying the prolongation of 
$\tZ$ to the generating
trivial cofibrations in $\spec(\simp)$. 
\end{proposition}

\begin{proof}
Define $\tU \mc \spec(\sab) \to \spec(\simp)$ by composing $U
\mc \spec(\sab) \to H\mZ\Mod$ 
with the forgetful functor from $H\mZ$-modules to underlying symmetric spectra.
The left adjoint of $\tU$ is defined by applying $\tZ$ to each level.  
Denote this prolongation by $\tZ$ as well.
%We apply the lifting theorem for cofibrantly generated model 
%categories from~\cite[2.3]{ss}.  
Define a new model structure on $\spec(\sab)$ with weak equivalences and 
fibrations the maps which are underlying weak equivalences and fibrations 
in $\spec(\simp)$.  The fibrant objects are then the $\Omega$-spectra and
the trivial fibrations are the level trivial
fibrations.  It follows that the cofibrations must be the
maps with the left lifting property with respect to the level trivial 
fibrations.  To see that these structures satisfy the axioms of a model
category we note that this structure agrees with   
the stable model structure defined in~\cite{hovey-stable}.  To see this,
note that the trivial fibrations agree.  Then the cofibrations and hence the 
cylinder objects must also agree.   Also, the fibrant objects agree and hence
the weak equivalences agree; see for example~\cite[7.8.6]{psh}. 

Since the weak equivalences and fibrations of $\spec(\sab)$ are
determined by $\tU$, they are detected by   
$\tZ(I)$ and
$\tZ(J)$ where $I=FI_{\partial}$ and $J=FI_{\Lambda} \cup K$ are 
the generating cofibrations and trivial cofibrations of $\spec(\simp)$
defined in~\cite[3.3.2, 3.4.9]{hss}.   
It follows that $\spec(\sab)$ is a cofibrantly generated
model category with generating cofibrations $\tZ(I)$ and 
generating trivial cofibrations $\tZ(J)$. 
%
%%%%%%%%%%%%%%%%5
%To apply this lifting theorem we verify
%the first criterion stated in~\cite[2.3]{ss}; we show that for
%the generating trivial cofibrations in $\spec$, pushouts and
%directed colimits of $\tZ(J)$ in $\spec(\sab)$ forget to weak equivalences in 
%$\spec$.  
%Both $\tZ$ and $\tU$ preserve level equivalences and level cofibrations, so
%pushouts and directed colimits of $\tZ(FI_{\Lambda})$ forget to
%level equivalences in $\spec$.  There is a map $H\mZ \sm_S F_{n+1}S^1 \to
%\tU\tZ F_{n+1}S^1$ which induces an isomorphism on stable homotopy
%groups and hence is a stable equivalence by~\cite[3.1.11]{hss}. 
%By the monoid axiom~\cite[5.4.1]{hss}, $H\mZ \sm_S -$ takes stably trivial 
%cofibrations
%to stable equivalences.   Then, by Ken Brown's lemma~\cite[1.1.12]{hovey-book},
%since $\lambda_n \mc F_{n+1}S^1 \to F_n S^0$ from~\cite[3.4.9]{hss} is
%a stable equivalence between stably cofibrant objects,  
%the map $H\mZ \sm_S F_{n+1}S^1 \to H\mZ \sm_S F_n S^0$ is also a  stable
%equivalence. 
%Since $\tU\tZ F_n S^0 \iso H\mZ \sm_S F_n S^0$, this implies that  
%$\tU\tZ\lambda_n$ is a stable equivalence for 
%$\lambda_n \mc F_{n+1}S^1 \to F_n S^0$. 
%It then follows that each map in $\tU\tZ(J)$
%is a level cofibration and a stable equivalence.  The necessary criterion
%follows since pushouts and directed colimits preserve such stable
%equivalences by the proof of~\cite[5.4.1]{hss}.  
\end{proof}

To show that the three adjoint pairs in Proposition~\ref{prop-list} 
are weak (or strong) monoidal Quillen equivalences, we use the following 
criterion for Quillen pairs between monoidal stable model categories 
from~\cite{emmc}.  Since
the unit objects are cofibrant in each case here, the unit condition
is simpler than the one appearing in~\cite{emmc}.   
An object $A$ in a stable model category $\C$ is said to
{\em stably detect weak equivalences} if $f\mc Y \to Z$ is a weak
equivalence if and only if $[A, Y]_k^{\HoC} \to [A, Z]_k^{\HoC}$ is
an isomorphism for all $k \in \mZ$.  

\begin{proposition}\label{prop-wm-crit}\cite[3.17]{emmc} 
Let $\lambda\mc \cD \adjoint \cC\mc R$ be a Quillen adjoint pair
between monoidal stable model categories  such that $R$ is lax
monoidal and the unit objects $\mIC$ and $\mID$ are both cofibrant.
Suppose further that
\begin{enumerate}
\item the adjoint to the monoidal structure map $\nu\mc \mID \to R(\mIC)$,
namely $\widetilde{\nu}\mc\lambda(\mID) \to \mIC$, is a weak equivalence in $\C$; and
\item the unit $\mID$ stably detects weak equivalences in $\D$. 
\end{enumerate}
Then $(\lambda, R)$ is a weak monoidal Quillen pair.
\end{proposition} 

Now we verify the three parts of Proposition~\ref{prop-list} separately
in Propositions~\ref{prop-hz-sab}, ~\ref{prop-sab-ch} and~\ref{prop-DR2}.

\begin{proposition}\label{prop-hz-sab}
%\item 
$(\spec(\sab), H\mZ\Mod, Z, U)$ is a strong monoidal Quillen equivalence.
Moreover,
$U$ preserves all weak equivalences. 
%so the fibrant condition in~\cite[3.12]{emmc} is unnecessary. 
\end{proposition}

\begin{proof}
%Since $U\Sym(\tZ S^1)\iso H\mZ$ and  
%the forgetful functor from $\sab$ to simplicial sets is lax symmetric
%monoidal, it follows that $U$ is lax symmetric monoidal.    
%
We first define the left adjoint $Z$.  As with $\tU$ in 
Proposition~\ref{prop-cof-gen}, one expects the left adjoint to involve $\tZ$. 
Applying $\tZ$ to each level of an $H\mZ$-module in $\spec(\sset)$ produces a 
$\tZ(H\mZ)$-module in $\spec(\sab)$.
The monad structure on $\tZ$ then induces a ring homomorphism 
$\mu\mc\tZ(H\mZ) \to H\mZ \iso \Sym(\tZ S^1)$ in $\spec(\sab)$
which induces a push-forward $\mu_*\mc \tZ(H\mZ)\Mod \to \spec(\sab)$.  
One can then check that $Z(X)= \mu_*(\tZ X) = \tZ X \sm_{\tZ(H\mZ)} H\mZ$ is
left adjoint to $U$.  Since both $\mu_*$ and $\tZ$ are strong symmetric 
monoidal so is $Z$.  It follows that $U$ is lax symmetric monoidal. 
Note that $\widetilde{\nu}\mc Z(H\mZ) \to \Sym(\tZ S^1)$ is an isomorphism and 
both units are cofibrant. 

It is clear that $U\mc \spec(\sab) \to H\mZ\Mod$ detects and preserves
weak equivalences and fibrations since in both cases these are determined
by the forgetful functor to $\spec(\simp)$; see Proposition~\ref{prop-cof-gen}. 
Hence, $U$ and $Z$ form a Quillen adjoint
pair.  Since $H\mZ$ stably detects weak equivalences in $H\mZ$-modules, it
follows by Proposition~\ref{prop-wm-crit} that $(U, Z)$ is a weak
monoidal Quillen pair.  Since $Z$ is strong monoidal, $(U,Z)$ is in
fact a strong monoidal Quillen pair.  

We are left with showing that this Quillen pair is a Quillen equivalence.
$U$ and $Z$ induce adjoint total derived functors $\widebar{U}$ and 
$\widebar{Z}$ on the homotopy categories.  
Since $U$ detects and preserves weak equivalences, by~\cite[4.1.7]{hss} we 
only need to check that $\psi \mc \id \to \widebar{U}\widebar{Z}$ is an 
isomorphism to establish the Quillen equivalence.  
First, note that $\psi$ is an isomorphism on $H\mZ$. 
These functors are exact and preserve coproducts, so
$\psi$ is an isomorphism on any object built from $H\mZ$ via suspensions,
triangles or coproducts; that is, $\psi$ is an isomorphism
on the localizing subcategory generated by
$H\mZ$.  Since $H\mZ$ detects the weak equivalences in $H\mZ$-modules,
it is a generator by ~\cite[2.2.1]{stable} and this localizing
subcategory is the whole category.
\end{proof}

\begin{proposition}\label{prop-sab-ch}
$( \spec(\sab), \spec(\ch), L, \phi^* N)$ is a weak monoidal Quillen equivalence.  
Moreover,
$\phi^* N$ preserves all weak equivalences. 
%so the fibrant conditions in~\cite[3.12]{emmc} are unnecessary. 
\end{proposition}

\begin{proof}
Since normalization from $\sab$ to $\ch$ is a lax
symmetric monoidal functor, its prolongation $N\mc \spec(\sab)
\to \N\Mod$ is also lax symmetric monoidal.  Recall, $\N= N(\Sym(\tZ S^1))
= N \tZ \bS$.
Since $\phi\mc\Sym(\mZ[1]) \to \N$ is a ring homomorphism between
commutative monoids, pulling back modules $\phi^*\mc
\N\Mod \to \spec(\ch)$ is also a lax symmetric
monoidal functor.  

We show that $\phi^*$ and  $N$ are each the right adjoint of a 
Quillen equivalence and each preserve all weak equivalences, 
thus $\phi^*N$ has these properties as well.  We delay the
verification that $(L,\phi^*N)$ is a weak monoidal Quillen pair until
the end of this proof.

Since in degree $n$ the map $\phi$ is 
induced by the monoidal structure on $N$, 
%$(\mZ[1])^{\otimes n} \iso 
$(N\tZ S^1)^{\otimes n} \to N(\tZ S^1{}^{\otimes n})$, 
it is a level weak equivalence by~\cite[29.4, 29.7]{may}.  
By the  Quillen invariance property, verified in 
Corollary~\ref{cor-mon-ax}, $\phi^*$
is then the right adjoint of a Quillen equivalence.
The model category on $\N\Mod$ has underlying
weak equivalences and fibrations, as defined by applying~\cite[4.1]{ss}.
Hence $\phi^*$ preserves all weak equivalences by definition.   

For the normalization functor, we first show that $N$ preserves all 
equivalences.  Since normalization takes weak equivalences in $\sab$ to
weak equivalences in $\ch$, $N$ preserves all levelwise weak equivalences.
Given any stable equivalence it can be factored as a stably trivial
cofibration followed by a stably trivial fibration.  Since stably
trivial fibrations are levelwise equivalences, we only need to show that
$N$ takes stably trivial cofibrations to stable equivalences.  We proceed
by showing that $N$ takes the generating stably trivial cofibrations
to stable equivalences which are also monomorphisms, that is trivial 
injective cofibrations.  Since $N$ commutes with colimits and trivial 
injective cofibrations are preserved under pushouts and directed colimits,
this then shows that $N$ takes each stably trivial cofibration to a
stable equivalence. 

Recall from Proposition~\ref{prop-cof-gen} that a set of generating stably 
trivial cofibrations for $\spec(\sab)$ is given by $\tZ(J)$ 
where $J= FI_{\Lambda} \cup K$ is defined in~\cite[3.4.9]{hss}.  
Here the maps in $FI_{\Lambda}$ are level equivalences and the maps
in $K$ are level equivalent to maps of the form 
$\lambda_n \sm X\mc F^{\simp}_{n+1}(X \sm S^1) \to F^{\simp}_nX$. 
Here we have added a superscript to denote $F_n^{\C}\mc \C \to \spec(\C)$
the left adjoint to evaluation at level $n$.  Since $N$ preserves 
level equivalences and monomorphisms, it is enough to show that $N$ of
the map $\tZ( \lambda_n \sm X)$ is a stable equivalence.
Note that in level $k$, $(N\tZ F^{\simp}_mY)_k \iso \Sigma_k 
\sm_{\Sigma_{k-m}} N \tZ(
Y \sm S^{k-m} ).$    The shuffle map induces maps 
$N\tZ X \otimes N \tZ A \to N\tZ(X \sm A)$ which are weak equivalences; 
see~\cite[2.7]{emmc}.   This in turn induces a levelwise weak equivalence 
 $\N \otimes_{\Sym(\mZ[1])} F^{\ch}_{n+1}(N\tZ X \otimes N\tZ S^1)
\to N \tZ F^{\simp}_{n+1}( X \sm S^1)$.  
Similarly, there is a map $\N \otimes_{\Sym(\mZ[1])} F^{\ch}_n(N\tZ X) \to
N \tZ F^{\simp}_n X$ which is also a levelwise weak equivalence.  
So the two horizontal maps in the diagram below are level weak equivalences.
\[
\begin{CD}
 \N \otimes_{\Sym(\mZ[1])} F^{\ch}_{n+1}(N\tZ X \otimes N\tZ S^1)
@>>> N \tZ F^{\simp}_{n+1}( X \sm S^1) \\ 
@VVV @VV{N(\lambda_n \sm X)}V \\
 \N \otimes_{\Sym(\mZ[1])} F^{\ch}_{n}(N\tZ X) 
@>>> N \tZ F^{\simp}_{n}( X )  
\end{CD}
\]
We next show the left vertical map is a stable equivalence.
Since $N \tZ X$
is cofibrant in $\ch$ and $N\tZ S^1\iso \mZ[1]$, \cite[8.8]{hovey-stable} 
shows that  $F^{\ch}_{n+1}(N\tZ X \otimes \mZ[1]) \to F^{\ch}_{n}(N\tZ X)$
is a stable equivalence.  
Since $\phi$ is a weak equivalence,
$\phi_* = \N \otimes_{\Sym(\mZ[1])} -$ is a left Quillen functor by 
Corollary~\ref{cor-mon-ax}. 
Hence, $\phi_*$ preserves stable equivalences between cofibrant objects
and the left vertical map is a stable equivalence.
We conclude that $N(\lambda_n \sm X)$ is level equivalent to a stable 
equivalence and hence is itself a stable equivalence.

We now show $N\mc \spec(\sab) \to
\N\Mod$ is the right adjoint of a Quillen equivalence.  
%$N$ is a right adjoint by Theorem~\ref{thm-aft} and Proposition~\ref{prop-loc};
Denote its left adjoint by $L'$.  Since $N \mc \sab \to \ch$ preserves  
weak equivalences and fibrations, the prolongation preserves
fibrations and weak equivalences between fibrant objects because they
are levelwise fibrations and levelwise weak equivalences~\cite[A.3]{dugger}.
By~\cite[A.2]{dugger}, this shows $N$ is a right Quillen functor.
Since weak equivalences between fibrant objects are level equivalences,
$N$ also detects such weak equivalences. 
So, by~\cite[4.1.7]{hss} we just
need to show that $\psi\mc X \to \widebar{N}\widebar{L'}X$ is an isomorphism
for all $\N$-modules $X$.  
First consider the unit $\N =N \tZ \bS$.  Since $N$ commutes with evaluation
at level zero (in $\sab$ or $\ch$), the associated left adjoints also
commute.   
\[\xymatrix@=10mm{ 
\spec(\sab) \quad \quad
\ar@<-.4ex>_-{\quad N\quad}[r] 
\ar@<-1ex>^-{ev_0}[d] 
&
\quad \quad \N\Mod 
\ar@<-.4ex>_-{L'}[l]
\ar@<3.8ex>^-{ev_0}[d]\\ 
\sab\quad \quad  
\ar@<1.8ex>^-{F_0^{\sab}}[u] 
\ar@<-.4ex>_-{N}[r] &
\quad \quad \quad  \ch 
\ar@<-.4ex>_-{\Gamma}[l]
\ar@<-3ex>^-{\N \otimes F_0^{\ch}(-)}[u] 
}
\]
So $L'$ of the free $\N$-module on $\mZ[0]$ is isomorphic
to the free $\Sym(\tZ S^1)$-module on $\Gam \mZ[0] = \mZ S^0$; that is, 
$L'\N=L'N\tZ \bS \iso \Sym(\tZ S^1) =\tZ \bS$.  See also~\cite[3.7]{emmc}.  
Since $N$ preserves all weak equivalences, 
$\widebar{N}\iso N$ and $\psi$ is an isomorphism on $N\tZ \bS$. 
Both model categories here are stable, so 
$\widebar{N}$ and $\widebar{L'}$ are both exact functors.  
Thus, $\psi$ is an isomorphism on the localizing subcategory generated
by $\N$. 
Since  $\N=N\tZ \bS$ is a generator for $\N$-modules, it follows
that $\psi$ is an isomorphism on all $\N$-modules. 

We use Proposition~\ref{prop-wm-crit} to establish 
$(\spec(\sab), \spec(\ch), \phi^* N)$ 
as a weak monoidal Quillen pair.  The unit objects
$\Sym(\bZ[1])$ and $\Sym(\tZ S^1)$ are both cofibrant and 
$\widetilde{\nu}\mc L'\phi_*
\Sym(\bZ [1]) \to \Sym(\tZ S^1)$ is an isomorphism since $\phi_*
\Sym(\bZ [1]) = \N$.   Finally, $\Sym(\bZ[1] ) = F_0^{\ch} \bZ [0]$ 
is a generator for $\spec(\ch)$ and thus detects weak equivalences.   
This follows from~\cite[2.2.1]{stable}
since $[F_0 \bZ[0], X]_k = [F_k \bZ[0], X] = [\bZ[0], X_k]$ and
weak equivalences between fibrant objects are level equivalences.    
\end{proof}

Next we define the left adjoint of the functor $R: \Ch \to \Spec(\ch)$.
Let $I$ be the skeleton of the category of finite sets and 
injections with objects ${\bf n}$.
Given $X$ in $\Spec(\ch)$ define a functor $\cD_X: I \to \Ch$ by
$\cD_X({\bf n})= \bZ[-n]\otimes X_n$.  For $X$ in $\spec(\ch)$
there is a structure map $\sigma\mc\bZ[m-n] \tensor X_n
\to X_m$ with adjoint $\widetilde{\sigma}\mc X_n \to \bZ[n-m] \otimes X_m$.  
For a standard inclusion of a
subset $\alpha: {\bf n} \to {\bf m}$ the map $\cD_X(\alpha)$ is 
$\bZ[-n]\otimes \widetilde{\sigma}.$ 
For an isomorphism in $I$, the action is given by the tensor
product of the action on $X_n$ and the sign action on $\bZ[-n]$.
The functor $D: \Spec(\ch) \to \Ch$ is defined by $DX = \colim_I \cD_X$.
This functor is similar to the detection functor introduced 
in~\cite[3.1.1]{b-thh}.
Now we show that $D$ is the left adjoint of $R$.

\begin{proposition}\label{prop-DR}
The functors $D$ and $R$ are adjoint.
\end{proposition}

To prove this we first calculate $D$ on free spectra.

\begin{lemma}\label{lem-D}
$D(F_mK) = \bZ[-m] \tensor K$. 
\end{lemma}

\begin{proof}
For $n \geq m$, 
$(F_mK)_n = \Sigma_n^+ \tensor_{\Sigma_{n-m}} \bZ[n-m] \tensor K$
and for $n < m$, $(F_mK)_n = 0$.  Since
$\hom_I({\bf m}, {\bf n}) \iso \Sigma_n/\Sigma_{n-m}$ as $\Sigma_n$ sets,
$(F_mK)_n$ can be rewritten as 
$\hom_I({\bf m},{\bf n}) \tensor \bZ[n-m] \tensor K$. 
Since $\hom_I({\bf m},-)$ is a free diagram (it is left adjoint to evaluation
of a diagram at ${\bf m}$), the proposition follows. 
\end{proof}

\begin{proof}[Proof of Proposition~\ref{prop-DR}]
Recall that $D^n$ is the acyclic chain complex modeling the $n$-disk. 
For $n > 0$, $\Hom_{\Spec(\ch)}(F_m D^n, RY)$ is the $n$th level of  $(RY)_m$,
and for $n = 0$ the zeroth level of  $(RY)_m$ is  
$\Hom_{\Spec(\ch)}(F_m \bZ[0], RY)$. 
Thus, it is enough to note that $\Hom_{\Ch}(D(F_m D^n), Y)
= \Hom_{\Ch}(D^{n-m}, Y) = Y_{n-m}$ and $\Hom_{\Ch}(D(F_m \bZ[0]), Y)
= \Hom_{\Ch}(\bZ[-m], Y) = Z^{-m}(Y)$ where $Z^{-m}(Y)$ is the kernel
of the differential $Y_{-m} \xrightarrow{d} Y_{-m-1}$.     
\end{proof}

\begin{proposition}\label{prop-DR2}
$(\Ch, \Spec(\ch), D, R)$ is a strong monoidal Quillen equivalence.  
Moreover, $R$ preserves all weak equivalences. 
\end{proposition}

\begin{proof}
First we show that $R$ preserves fibrations.  Actually, we show that $R$ takes 
fibrations to level fibrations between fibrant objects.  Since the
stable model category on $\Spec(\ch)$ is a localization of the levelwise
structure, these are in fact fibrations by~\cite[A.3]{dugger}.  It is easy to
check that $R$ takes fibrations in $\Ch$ to level fibrations in  
$\Spec(\ch)$ because a fibration in $\Ch$ is a surjective map and a fibration 
in $\ch$ is a surjection above degree zero.  
Recall that a fibrant object $X$ in $\Spec(\ch)$ is the analogue of an 
$\Omega$-spectrum; each $X_n$ is levelwise fibrant and $X_n \to C_0(\bZ[-1] \tensor
X_{n+1})$ is a quasi-isomorphism.  Since all objects in $\ch$ are fibrant,
the only condition to verify follows from checking that $(RY)_n \iso
C_0(\bZ[-1] \tensor (RY)_{n+1})$. 

Next we show that $R$ preserves and detects weak equivalences.  If
$X \to Y$ is a quasi-isomorphism in $\Ch$, then $C_0(X \tensor \bZ[n])
\to C_0(Y \tensor \bZ[n])$ is also a quasi-isomorphism for any $n$.
Thus $RX \to RY$ is a level equivalence and hence also a stable equivalence. 
Also, if $RX \to RY$ is a stable equivalence, then it must be
a level equivalence since both $RX$ and $RY$ are fibrant.  Since 
$H_n(RX)_m = H_{n-m}X$ for $n \geq 0$, it follows that $X \to Y$ is
a quasi-isomorphism.

Since $R$ preserves and detects weak equivalences, to show that
$D$ and $R$ form a Quillen equivalence it is enough to check that the derived
adjunction is an isomorphism on a generator as in the proofs of
Proposition~\ref{prop-hz-sab} and~\ref{prop-sab-ch}.  
Since $F_0({\bZ}[0])$ is
a generator for $\Spec(\ch)$ and   $RD(F_0K)=RK = F_0 K$ for $K$ in $\ch$, 
this follows.

Next we show that $D$ is strong symmetric monoidal; it then follows that 
$R$ is lax symmetric monoidal.  First, $D$ is symmetric monoidal because 
there is a natural transformation 
$\gamma_{X,Y}: DX \otimes DY \to D(X \sm Y)$.  Since colimits
commute with tensor, the source can be rewritten as $\colim_{I \times I}
(X_n \otimes Y_m \otimes \bZ[-n-m])$.  The colimit
in the target can be pulled back over the functor $p:I\times I \to I$
given by $p(n,m) = n +m$.  Then $\gamma_{X,Y}$ is induced by the
inclusion $X_n \otimes Y_m \to (X \sm Y)_{n+m}$.
To show that $\gamma$ is always an isomorphism we first verify
this for free spectra.  By Lemma~\ref{lem-D}, 
$D(F_mK) \otimes D(F_nL) = \bZ[-m-n]\otimes K \otimes L$ which is isomorphic
to $D(F_mK \sm F_nL) = DF_{n+m}(K\otimes L)$.  It follows that
$\gamma$ is always an isomorphism since any spectrum 
$Z$ is the coequalizer of the two maps from $FFZ$ to $FZ$ where 
$FZ = \oplus_n F_n(Z_n)$.    
Since the unit
objects $F_0 \bZ[0]$ and $\bZ[0]$ are both cofibrant, $D(F_0\bZ[0])=
\bZ[0]$ and $\bZ[0]$ stably detects weak equivalences,
it follows from Proposition~\ref{prop-wm-crit} that this is a 
strong monoidal Quillen equivalence.  
\end{proof}

We now turn to the proof of Theorem~\ref{thm-Q-comm} which states a partial 
extension 
of the main result of this paper to the commutative case over the rationals.
We first need the following fact about $D$. 

\begin{lemma}\label{lem-D-rat}
$D$ preserves all weak equivalences over $\bQ$.   
\end{lemma}

\begin{proof}
Any weak equivalence can be factored as a trivial cofibration
followed by a trivial fibration.  $D$ preserves trivial cofibrations,
since $R$ preserves fibrations.  Any trivial fibration is in fact
a level equivalence, so to show that $D$ takes trivial fibrations
to weak equivalences we only need to show that rationally 
$DX \iso \hocolim_I( X_n \tensor \bZ[-n])$.   

First, note that $\hocolim_I = \colim_n (\hocolim_{I_n})$ where 
$I_n$ is the full subcategory with objects ${\bf i}$ for $i \leq n$. 
Similarly $\colim_I = \colim_n(\colim_{I_n})$, so we just need
to consider $\hocolim_{I_n}$.  Consider the category $I/n$ of
objects over ${\bf n}$.  There is a $\Sigma_n$ action on
$\hocolim_{I/n}$ coming from the action on ${\bf n}$.  For any diagram
$F$ the quotient $(\hocolim_{I/n} F)/ \Sigma_n$ is isomorphic to
$\hocolim_{I_n} F$.  Also $\hocolim_{I/n} F \to F({\bf n})$ is a
$\Sigma_n$ equivariant map and a quasi-isomorphism since $\id:{\bf n}
\to {\bf n}$ is the final object in $I/n$.  Since taking the quotient by 
this $\Sigma_n$ action is exact rationally, 
$(\hocolim_{I/n} F) /\Sigma_n \to F({\bf n})/\Sigma_n$ is a quasi-isomorphism
over $\bQ$.  Thus, $\hocolim_{I_n} F \to F({\bf n})/\Sigma_n = \colim_{I_n} F$
is a quasi-isomorphism.
\end{proof}

\begin{proof}[Proof of Theorem~\ref{thm-Q-comm}]
Throughout this proof let $C$ be a commutative $H\bQ$-algebra spectrum.
We consider the restriction of $\pPsi =Dc \phi^*N Zc$ to commutative 
$H\bQ$-algebra spectra and show that for each such $C$, $\pPsi C$ is naturally 
weakly equivalent to a commutative differential graded $\bQ$-algebra. 
The three main functors $D$, $\phi^*N$ and $Z$ in $\pPsi$ are lax symmetric 
monoidal functors and hence strictly preserve commutative rings.  The cofibrant
replacement functors of monoids involved in $\pPsi$ are not symmetric monoidal 
though.  This is why $\pPsi C$ is only weakly equivalent and not isomorphic
to a commutative dg $\bQ$ algebra.   
 
Since $D$ rationally preserves all weak equivalences by 
Proposition~\ref{lem-D-rat}, the cofibrant replacement transformation
$c \to \id$ induces a natural weak equivalence
$\pPsi  \to \widetilde{\pPsi}=  D \phi^* N Zc$  
when restricted to $H\bQ$-algebras.  To complete this proof we show that the 
functor $\widetilde{\pPsi}$ is related by a zig-zag of natural
transformations, each of which induces a weak equivalence on any 
$H\bQ$-algebra, to a lax symmetric monoidal functor $\pPsi'$. 
%Hence, $\pPsi C $ is naturally weakly
%equivalent to $\pPsi' C$ which is a commutative dg $\bQ$-algebra. 

We have one remaining cofibrant replacement functor to consider in 
$\widetilde{\pPsi}$.   
We show that there is a zig-zag of natural transformations between $Z c$
and a lax symmetric monoidal functor $\alpha^* \bQt$ which induces
weak equivalences on $H\bQ$-algebra spectra.  
Define $\pPsi' = D  \phi^* N \alpha^* \bQt$; since each factor
is lax symmetric monoidal so is $\pPsi'$. 
Since $D$ and $\phi^*N$ rationally preserve all weak equivalences,  
the zig-zag between $Z c$ and $\alpha^* \bQt$ induces a zig-zag
of weak equivalences between $\pPsi'C$ and $\widetilde{\pPsi} C$.  Hence 
$\pPsi C$ is also naturally weakly equivalent to $\pPsi'C$, a rational 
commutative dga, for $C$ any $H\bQ$-algebra spectrum.

As above let $\tZ\mc \simp \to
\sab$ be the free abelian group functor 
on the non-basepoint simplices.   Define $\bQt$ similarly.
We must compare $Z= \bZt(-) \sm_{\bZt(H\bZ)} H\bZ$, defined in detail in
the proof of Proposition~\ref{prop-hz-sab}, 
with its rational analogues.  Define $\widebar{Q}=
\bQt(-)\sm_{\bQt(H\bZ)} H\bQ$ and $Q= \bQt(-) \sm_{\bQt(H\bQ)} H\bQ$.
Then $\widebar{Q} = Q \circ [(-) \sm_{H\bZ} H\bQ]$ because $\bQt$ is
strong monoidal.
The inclusion $\bZ \to \bQ$ induces a natural monoidal transformation 
$Z \to \widebar{Q}.$  One can check that
for $c$ a cofibrant replacement functor for $H\bZ$-algebras,
$c_{\bQ}= c(-) \sm_{H\bZ} H\bQ$ is a cofibrant replacement
functor of $H\bQ$-algebras.  This follows since $c \to c_{\bQ}$ is a weak
equivalence on $H\bQ$-algebras by~\cite[IV.4.1]{ekmm} for example,
since $\bQ$ is flat over $\bZ$.  So there is a natural transformation 
$Z c \to \widebar{Q} c \xrightarrow{\iso} Q c_{\bQ}$.  
Proposition~\ref{prop-hz-sab} shows
that $Z$ and $U$ form a Quillen equivalence.  Thus, since $U$ preserves all
weak equivalences, $c \to U Z c$ induces a weak equivalence on any
$H\bZ$-algebra.  Similarly 
one can show that $c_{\bQ} \to U Q c_{\bQ}$ induces a weak equivalence on
any $H\bQ$-algebra.  Since $U$ detects weak equivalences it follows  
that $Zc A \to Qc_{\bQ}A$ is a natural weak equivalence of algebras on any 
$H\bQ$-algebra $A$.
    
Now consider the maps $H\bQ \xrightarrow{\alpha} \bQt(H\bQ)
\xrightarrow{\beta} H\bQ$ given by the unit and multiplication of the
monad structure on $\bQt$.  
Both of these maps induce isomorphisms on
$\pi_*$, so they induce two Quillen equivalences via extension and
restriction of scalars, 
($\alpha_*$, $\alpha^*$), ($\beta_*$, $\beta^*$), between the respective
categories of algebras over $\spec(\sab)$~\cite[5.4.5]{hss}.  
Since $\beta \circ \alpha = \id$, $\alpha^* \beta^* =  \id$. 
So the functor $Q$ can be rewritten as 
$Q \iso \beta_* \bQt \iso \alpha^*\beta^*\beta_*\bQt$.  
The functor $\bQt\mc \spec(\simp) \to \spec(s\bQ\Mod)$ is a left Quillen 
functor since its right adjoint,
the forgetful functor, preserves weak equivalences and fibrations. 
So $\bQt$ preserves cofibrant objects and the Quillen equivalence
($\beta_*$, $\beta^*$) induces a weak equivalence $\bQt c_{\bQ}
\to \beta^*\beta_*\bQt c_{\bQ}$.  Since $\alpha^*$ preserves all
weak equivalences, this gives a natural weak equivalence $\alpha^* \bQt c_{\bQ}
\to \alpha^* \beta^* \beta_* \bQt c_{\bQ} \iso Q c_{\bQ}$. 
Finally, since $\bQt$ also preserves all weak equivalences,
there is a natural weak equivalence $\alpha^* \bQt c_{\bQ} \to \alpha^* \bQt$.
This produces the promised zig-zag between $Zc$ and $\alpha^* \bQt$.
Notice, since $\alpha^*$ is lax symmetric monoidal and $\bQt$ is strong 
symmetric monoidal functor, $\alpha^* \bQt$ is lax symmetric monoidal 
as required.  

%The required functor $f'$ is constructed as the fibrant replacement
%functor for a model category on the commutative monoids in 
%$\spec(\chZ_{\bQ})$ where a map is a fibration or weak equivalence if
%it is such on the underlying model category.  Then for any commutative
%monoid $C$, since $f'C$ is also fibrant as an underlying monoid, the 
%lifting property implies one can extend the weak equivalence
%$C \to f'C$ over the trivial cofibration $C \to fC$.  This gives a weak
%equivalence $fC \to f'C$.
%Since this is a weak equivalence between underlying fibrant objects,
%it is a level equivalence and hence $\Ev f C \to \Ev f' C$ is also a
%weak equivalence.  Note, similar arguments do not work for cofibrant
%replacements because a cofibrant commutative monoid need not be cofibrant
%as an underlying monoid. 
%
%The arguments for establishing this model category on the commutative
%monoids in $\spec(\chZ_{\bQ})$ follow just as for the commutative
%symmetric ring spectra in~\cite[Section 15]{mmss}.  Here though, one
%does not have to restrict to the positive model category because in
%$\chZ_{\bQ}$ the map $E\Sigma_i \to \bQ$ is a $\Sigma_i$-weak equivalence.
%Hence the key lemma,~\cite[15.5]{mmss}, applies for $n \geq 0$ rather
%than $n > 0$ and for any cofibrant object in $\spec(\chZ_{\bQ})$. Thus
%the fibrations and weak equivalences can be defined as in the ordinary
%underlying model category for $\spec(\chZ_{\bQ})$ rather than
%the positive model category.
\end{proof}

Finally, we consider the two step alternative comparison of 
$\spec(\ch)$ and $\chZ$.  

\begin{proposition}\label{prop-extra}
$(\spec(\chZ), \spec(\ch), i, \co)$ and $ (\spec(\chZ),  \chZ, F_0, \Ev)$ 
are strong monoidal Quillen equivalences. 
\end{proposition}

\begin{proof}
Since fibrations in $\ch$ are the maps which are surjections above
degree zero, $\co \mc \chZ \to \ch$ preserves fibrations and weak equivalences. 
Thus, $\co$ preserves fibrations and weak equivalences between stably
fibrant spectra because they are levelwise fibrations and levelwise weak 
equivalences by~\cite[A.3]{dugger}.
So by the criterion given in~\cite[A.2]{dugger}, $i$ and $\co$ form
a Quillen adjunction on the stable model categories.  
For an $\Omega$-spectrum $X$ in $\spec(\Ch)$ each negative homology group at 
one level is isomorphic to a non-negative homology group at a higher level,
$H_{-k} X^n \iso H_0 X^{n+k}$ for $k \geq 0$.  
Thus the functor $\co$ also detects weak equivalences between fibrant objects.  
By~\cite[4.1.7]{hss}, to show ($i, \co$) is a Quillen equivalence, 
it is thus enough to check 
that the derived adjunction is an isomorphism on the generator.  Since
$\Sym(\mZ[1])$ is concentrated in non-negative degrees, this follows. 

Since the inclusion $i\mc \ch\to \chZ$ is strong
symmetric monoidal and $i(\Sym(\mZ[1]) \iso \Sym (\mZ[1])$, the prolongation 
of $i$ is also strong symmetric monoidal.  Similarly, $\co \mc \chZ \to \ch$
is lax symmetric monoidal and so is its prolongation.  Since 
%$\Sym(\mZ[1])$
%is a generator of $\spec(\ch)$ by Proposition~\ref{prop-sab-ch}, 
%$\eta\mc\co(Sym(\mZ[1])
%\xrightarrow{\iso} Sym(\mZ[1])$ and 
both units are cofibrant, the first pair is a strong monoidal Quillen pair.

Since $\chZ$ is a stable model category, the second pair of adjoint functors 
form a Quillen equivalence by~\cite[9.1]{hovey-stable}.  
Both $F_0$ and $\Ev$ are strong symmetric monoidal functors.
Since 
%$[\mZ[0],X]_* = H_*X$, $\mZ[0]$ is a generator for
%$\chZ$ by~\cite[2.2.1]{stable}.  Also $\eta\mc
%F_0\mZ[0] \xrightarrow{\iso} \Sym(\mZ[1])$ and 
both units are cofibrant, 
the second pair is also a strong monoidal Quillen pair.
\end{proof}

\end{document}